%% file: main.tex
\documentclass[11pt]{amsart}
\usepackage[foot]{amsaddr}
\input{preamble.tex}

\title[HDA I]{Higher-Dimensional Algebra I: \\
       Braided Monoidal 2-Categories}

\author[Baez]{John C.\ Baez$^1$} 
\address{$^1$Department of Mathematics, University of California, Riverside CA, 92521, USA}

\author[Neuchl]{Martin Neuchl$^2$}
\address{$^2$Mathematisches Institut der Universit\"at M\"unchen,
Theresienstr.\  39,  80333 M\"unchen,
Germany}

\email{baez@math.ucr.edu, neuchl@rz.mathematik.uni-muenchen.de}

\date{October 2, 1995}

\begin{document}

\begin{abstract}
    We begin with a brief sketch of what is known and conjectured concerning braided monoidal 2-categories and their relevance to 4d TQFTs and 2-tangles.   Then we give concise definitions of semistrict monoidal 2-categories and braided monoidal 2-categories, and show how these may be unpacked to give long explicit definitions similar to,  but not quite the same as, those given by Kapranov and Voevodsky. Finally, we describe how to construct a semistrict braided monoidal 2-category $\Z(\C)$ as the `center' of a semistrict monoidal category $\C$, in a manner analogous to the construction of a braided monoidal category as the center of a monoidal category.  As a corollary this yields a strictification theorem for braided monoidal 2-categories.
\end{abstract}

\maketitle

\input{1intro.tex}

\input{2definitions.tex}

\input{3centerconstruction.tex}

\input{4embedding.tex}

\input{5conclusions.tex}

\subsection*{Acknowledgements}
Many of the basic ideas behind this paper were developed in collaboration with James Dolan.  We would also like to thank Lawrence Breen, Misha Kapranov, John Power, James Stasheff and Ross Street for useful discussions and correspondence.  J.B.\ is grateful to the Mathematical Institute of the University of Munich for their 
hospitality while some of this work was done.  We thank Joe Moeller and especially Jason Erbele for typesetting this new version of the paper in 2019.

\end{document}

%% file: preamble.tex
\usepackage[utf8]{inputenc}
\usepackage{color}
\usepackage{amsmath,amsthm}
\usepackage{amsfonts,amssymb,mathrsfs}
\usepackage[a4paper,top=3cm,bottom=3cm,inner=3cm,outer=3cm]{geometry}
\usepackage{tikz}
\usepackage{tikz-cd}

\renewcommand{\hom}{{\rm hom}}
\newcommand{\gtimes}{\otimes_{\rm G}}
\newcommand{\op}{{\rm op}}
\newcommand{\id}{{\rm id}}
\newcommand{\maps}{\colon}
\newcommand{\define}[1]{{\bf \boldmath{#1}}} 

\newcommand{\C}{\mathcal C}
\newcommand{\D}{\mathcal D}
\newcommand{\F}{\mathcal F}
\newcommand{\G}{\mathcal G}
\newcommand{\I}{\mathcal I}
\newcommand{\Z}{\mathcal Z}

\newcommand{\namedcat}[1]{\mathsf{#1}}
\newcommand{\Cat}{\namedcat{Cat}}
\newcommand{\Cob}{\namedcat{Cob}}
\newcommand{\Reps}{\namedcat{Reps}}
\newcommand{\Set}{\namedcat{Set}}
\newcommand{\Vect}{\namedcat{Vect}}

\newcommand{\stri}{-}
\newcommand{\tilR}{\tilde R}
\newcommand{\tilT}{\tilde T}
\newcommand{\zwei}{\Downarrow}
\newcommand{\To}{\Rightarrow}
\newcommand{\iso}{\cong}

\theoremstyle{plain}
\newtheorem{thm}{Theorem}

\newtheorem{lma}[thm]{Lemma}
\newtheorem{bem}[thm]{Remark}

\theoremstyle{definition}
\newtheorem{defn}[thm]{Definition}

\definecolor{myurlcolor}{rgb}{0.6,0,0}
\definecolor{mycitecolor}{rgb}{0,0,0.8}
\definecolor{myrefcolor}{rgb}{0,0,0.8}
\usepackage[pagebackref]{hyperref}
\hypersetup{colorlinks, linkcolor=myrefcolor, citecolor=mycitecolor, urlcolor=myurlcolor}

%% file: 1intro.tex
\section{Introduction}

This is the first of a series of articles developing the program introduced in the paper `Higher-Dimensional Algebra and Topological Quantum Field Theory' \cite{BD}, henceforth referred to as `HDA'.   This program consists of generalizing algebraic concepts from the context of set theory to the context of $n$-category theory, and using the resulting language to unify topological quantum field theory with traditional algebraic topology.  Rather than doing so systematically from the ground up, the papers in this series will instead address specific issues as they become manageable.  The present paper treats a concept which appears to be of special interest in  4-dimensional topology and physics: that of a braided monoidal 2-category.

To understand this concept and its role in higher-dimensional algebra, it is useful to recall some ideas described more thoroughly in HDA.  Loosely speaking, an $n$-category is a structure generalizing a category in which there are 0-morphisms or `objects', 1-morphisms between objects, 2-morphisms between 1-morphisms, and so on up to $n$-morphisms.  Giving a precise and sufficiently general definition of $n$-categories is, however, a rather subtle matter.  So-called `strict' $n$-categories can already be defined recursively for all $n$, using the idea that for any two objects $A$ and $B$ of an $n$-category,  $\hom(A,B)$ should be not a set but an $(n-1)$-category.   One can also unpack this recursive definition and obtain a definition in terms of an explicit list of operations for composing $j$-morphisms and  equational laws the operations obey \cite{Street}.   

However, strict $n$-categories violate the fundamental principle that \emph{``In any category it is unnatural and undesirable to speak about equality of two objects''} \cite{KV}.  It is all too easy to mistakenly treat two objects of a category as `equal' when they are merely
isomorphic, so it is better to systematically avoid such mistakes
by replacing all equations by specified isomorphisms.  Of
course, an isomorphism satisfies equations of
its own, which state that it is invertible, and in a 2-category these
equations themselves should be replaced by specified 2-isomorphisms, and
so on.  This leads to the recursively defined
notion of an `equivalence': a $j$-morphism that is 
strictly invertible if $j = n$, but only invertible up to an 
equivalence if $j < n$.   The practical advantages of replacing equations
by specified equivalences are already quite apparent in homotopy
theory, and they are likely to
become increasingly evident in other branches of mathematics and physics,
such as topological quantum field theory.
 
Taking this philosophy seriously, it is clear that
one should define a notion of `weak' $n$-category by
taking the definition of strict $n$-category and replacing all
equational laws between $j$-morphisms (for $j < n$) by specified
equivalences.   To serve essentially the same role as the
equations they replace, these equivalences should satisfy
some `coherence laws'.  However, to follow the weakening
principle, these `laws' should themselves not
be equations, in general, but only specified equivalences, and so
on: true equational laws are only to be required at the level of
$n$-morphisms.  Unfortunately, determining the correct coherence laws is 
a rather tricky business, so that weak $n$-categories have been defined so
far only for $n \le 3$.   They are usually called bicategories
\cite{Benabou} for $n = 2$ and tricategories \cite{GPS}
for $n = 3$.  A major challenge for higher-dimensional algebra
is to find a good theory of weak $n$-categories for all $n$.

In any event, one expects quite generally that in either the
strict or the weak context an $(n+1)$-category $\tilde \C$ with only one object
$\ast$ can be regarded as an $n$-category $\C$ by re-indexing, the
$j$-morphisms of $\C$ being simply the $(j+1)$-morphisms of
$\tilde \C$.  The $n$-categories we obtain this way have extra
structure.   For example, since the objects of $\C$ are
really morphisms in $\tilde \C$ from $\ast$ to itself, we can
`tensor' or compose them.  A category equipped with tensor
products is known as a monoidal category, and by analogy we
call any $n$-category arising from an $(n+1)$-category with one
object in this way a `monoidal $n$-category'.   Strict monoidal $n$-categories
are well understood for all $n$, while the weak ones are presently
defined only for $n \le 2$, since weak $n$-categories are only
defined for $n \le 3$.  

Similarly, we expect that an $(n+2)$-category $\tilde \C$ 
with only one object $\ast$ and one 1-morphism $1_\ast$
can be regarded as an $n$-category $\C$ with still further structure.
In particular, the tensor product should satisfy a kind of
commutativity condition.  When $n = 0$, this commutativity
condition is simply the equation $x \otimes y = y \otimes x$,
and it follows from a beautiful argument used by
Eckmann and Hilton \cite{EH} to show the commutativity of the higher
homotopy groups.   For simplicity, let $\C$ be a 
strict 2-category with
only one object $\ast$ and one 1-morphism $1_\ast$.  Then for any 
2-morphisms $x$ and $y$ in $\tilde \C$, both the horizontal composite
$x \otimes y$ and the vertical composite $xy$ are well-defined.
The 2-morphism $1 = 1_{1_\ast}$ is the unit for both vertical and
horizontal composition, and the exchange identity 
\[
(x_1 \otimes x_2)(y_1 \otimes y_2) = (x_1 y_1) \otimes (x_2y_2) 
 \]
holds for all 2-morphisms $x_i, y_i$.  Thus we have
\begin{align*}   x \otimes y & = (x1) \otimes (1y) \\
    & = (x \otimes 1)(1 \otimes y) \\
    & = xy \\
    & = (1 \otimes x)(y \otimes 1) \\
    & = (1y) \otimes (x1) \\
    & = y \otimes x,
\end{align*}
so vertical and horizontal composition are equal and
$\C$ is a commutative monoid.  Conversely, 
any commutative monoid is the set of
2-morphisms in some 2-category with one object and one 1-morphism.

As a consequence of the philosophy underlying weak $n$-categories,
when $n = 1$ this commutativity condition is not an equation but
an isomorphism.  In other words, a weak 3-category with only
one object and one 1-morphism can be thought of as a weak `braided'
monoidal category: one equipped with a natural isomorphism
\[           
    R_{x,y} \maps x \otimes y \to y \otimes x 
\]
satisfying certain coherence laws \cite{GPS,JS}.  More generally,
we may define a `braided monoidal $n$-category' to be an $(n+2)$-category
with one object and one 1-morphism.   More generally still, a 
 $(n+k)$-category with only one $j$-morphism
for each $j < k$ can be regarded as a special sort of $n$-category,
a `$k$-tuply monoidal $n$-category'.  These
play a key role in HDA, from which the table in Figure 1 is taken.
Note in particular the `stabilization' predicted for $k \ge n+2$.

\begin{center}
\begin{figure}[ht]
{\small
\begin{tabular}{|c|c|c|c|}  \hline
      & $n = 0$   & $n = 1$    & $n = 2$          \\     \hline
$k = 0$  & sets      & categories & 2-categories     \\     \hline
$k = 1$  & monoids   & monoidal   & monoidal         \\    
      &           & categories & 2-categories     \\     \hline
$k = 2$  &commutative& braided    & braided          \\      
      & monoids   & monoidal   & monoidal         \\
      &           & categories & 2-categories     \\     \hline
$k = 3$  &`'         & symmetric  & weakly involutory \\   
      &           & monoidal   & monoidal         \\
      &           & categories & 2-categories     \\     \hline
$k = 4$  &`'         & `'         &strongly involutory\\      
      &           &            & monoidal         \\
      &           &            & 2-categories     \\     \hline
$k = 5$  &`'         &`'          & `'               \\
      &           &            &                  \\ 
      &           &            &                  \\      \hline
\end{tabular}} \vskip 1em
\caption{Weak $k$-tuply monoidal $n$-categories: expected results}
\end{figure}
\end{center} 

Unfortunately, the weak versions of these structures have only
been defined in certain cases so far.  In particular, the weak
version of braided monoidal 2-categories is not yet  understood,
because they should be weak 4-categories with only one object and
one 1-morphism, and weak 4-categories have not yet been defined.
However, Kapranov and Voevodsky \cite{KV} have defined a more limited
class of `semistrict' braided monoidal 2-categories, 
the hope being that eventually all weak braided monoidal 2-categories
could be proven equivalent to these semistrict
ones (in some appropriate sense).   This strategy has already proven
successful at other levels.   For example, Gordon, Power, and Street \cite{GPS}
showed that all weak 3-categories are equivalent to a certain
class of semistrict ones, and as a corollary, 
all weak monoidal 2-categories are equivalent to certain semistrict ones.  
Since braided monoidal 2-categories can
be thought of as monoidal 2-categories equipped with extra structure, 
one expects a similar `strictification
theorem' to hold at the level of braided monoidal 2-categories.

Kapranov and Voevodsky's definition of a semistrict braided monoidal
2-category consists of a long explicit list of operations and
equational laws.   The first main goal of this paper is to
present a more concise and conceptual definition.  When
we unpack this definition to obtain an explicit list of
operations and laws, we find that it differs from Kapranov and
Voevodsky's list in a few places.  These appear to be slight defects
in their definition; for example, our subsequent theorems would
not work as smoothly if we used their definition.

\subsection{The Center Construction}
\label{center}

The second main goal of this paper is to give a procedure for
constructing a braided monoidal 2-category as the `center' $\Z(\C)$ of
of a monoidal 2-category $\C$.  To appreciate this
rather complicated procedure it is necessary to
understand the general concept of `center' proposed in HDA.
In essence this concept is simple; all the complications
arise from the lack of a good general theory of weak $n$-categories.  

There is no `set of all sets', but there is a class of all sets.
Better still, there is a category $\Set$ having sets as
objects and functions between them as morphisms.  Similarly,
there is a 2-category $\Cat$ having small categories as
objects,  functors between them as 1-morphisms, and 
natural transformations between functors as 2-morphisms.
In general, we expect there
to be a very important $(n+1)$-category $n\Cat$ having as objects
all small $n$-categories  (i.e., those for which the $j$-morphisms
form a set).    This has been worked out quite
generally in the strict context, but in the weak context only
for $n \le 2$ \cite{Benabou,GPS}.

In terms of this idea, the `center' of a small $k$-tuply monoidal
$n$-category $\C$ is a small $(k+1)$-tuply monoidal $n$-category
$\Z(\C)$ defined as follows.  Recall that $\C$ is really a
special sort of $(n+k)$-category, namely one with only one
$j$-morphism for $j < k$.  Thus $\C$ is an object in $(n+k)\Cat$.  
Let $1_1 = 1_\C$ denote the identity  1-morphism of
$\C$ in $(n+k)\Cat$, and recursively define
\[   1_{j+1} = 1_{1_j} ,\]
so that $1_j$ is a $j$-morphism.   Then there should be a
sub-$(n+k)$-category $\Z(\C)$  of $(n+k)\Cat$ having $\C$ as its
only object, $1_\C$ as its only 1-morphism, $1_{1_\C}$ as its only
2-morphism, and so on up to $1_k$, and then having all
$(k+1)$-morphisms from $1_k$ to itself as $(k+1)$-morphisms, all
$(k+2)$-morphisms between these as $(k+2)$-morphisms, and so on.  
Since $\Z(\C)$ has only one $j$-morphism for $j < k+1$, it follows that
$\Z(\C)$ is a $(k+1)$-tuply monoidal $n$-category.  

As this construction is a bit mind-boggling at first sight, let
us illustrate it in the case $n = 0$, $k = 1$.  Thus we begin
with a small category $\C$ with only one object $\ast$.   
The set $\tilde \C$ of 1-morphisms of $\C$ can be an arbitrary monoid.  
Similarly, $\Z(\C)$ is a 2-category with only one object and
one 1-morphism, and the 2-morphisms of such a 2-category form a
commutative monoid.   
More precisely, $\Z(\C)$ is the sub-2-category of $\Cat$ having
$\C$ as its only object, $1_\C$ as its only 1-morphism, and
all natural transformations $T \maps 1_\C \to 1_\C$ as 2-morphisms.
What is such a natural transformation in concrete terms?  It
must assign to the one object $\ast$ of $\C$
a morphism $T_\ast \maps \ast \to \ast$, such that for 
all $f \maps \ast \to \ast$ the following diagram commutes:

\[
\begin{tikzcd}
    \ast
    \arrow[r, "f"]
    \arrow[d, "T_\ast", swap]
    &
    \ast
    \arrow[d, "T_\ast"]
    \\
    \ast
    \arrow[r, "f", swap]
    &
    \ast
\end{tikzcd}\]

In other words, it is simply an element $T_\ast$ of the 
center of $\tilde \C$.   Thus the generalized concept of center
reduces in this case to the standard notion. 

The case $n = 1$, $k = 1$ is more interesting.  The center of a weak monoidal category is a weak braided monoidal category \cite{JS,KV,Majid}.  In particular, if $H$ is a Hopf algebra, the category $\Reps(H)$  of finite-dimensional comodules of $H$ is a weak monoidal category, and the center $\Z(\Reps(H))$ is then the category of representations of a coquasitriangular Hopf algebra $DH$ called the `quantum double' of $H$.   (Working with comodules and coquasitriangular Hopf algebras, rather than modules and quasitriangular Hopf algebras, serves as a technical convenience.)  The
quantum double construction, invented by Drinfeld \cite{Drinfeld},
gives to many interesting coquasitriangular Hopf algebras.
In particular, the quantum groups arising from semisimple
Lie groups, while not quantum doubles themselves,
are straightforward quotients thereof \cite{Kassel}.  Thus the
center construction can be regarded as an elegant
approach to quantum groups, which, as we shall see, 
makes their appearance in 3-dimensional topology much less mysterious.   

The class of theorems known as `Tannaka--Krein reconstruction
theorems' \cite{DM,Majid2,Ulbrich} further
clarifies the relation between the center construction and 
quantum doubles.  Given a Hopf algebra $H$, the category
$\Reps(H)$ is a $\mathbb{C}$-linear abelian rigid monoidal category
and equipped with a faithful $\mathbb{C}$-linear exact monoidal functor to $\Vect$.
Conversely, given any such category $\C$ equipped
with such a functor to $\Vect$, $\C$ is equivalent to $\Reps(H)$
for some Hopf algebra $H$ unique up to natural isomorphism.  
A similar theorem holds for $H$ coquasitriangular and $\C$ braided.  
Thus we may construct the quantum double of $H$ by first forming
$\Reps(H)$, then taking the center $\Z(\Reps(H))$ of this category, and then
applying Tannaka--Krein reconstruction to obtain $DH$.

It is natural to hope that other cases of the center
construction will give interesting analogs of these results.
The most interesting
case that can be handled with our present limited understanding
of weak $n$-categories is the case $n = 2$, $k = 1$:
if $\C$ is a monoidal 2-category, one expects that
$\Z(\C)$ will be a braided monoidal 2-category.   The difficulty with
proving this result is that we lack a general theory of weak
4-categories.  Thus we do not know the definition of a weak braided monoidal
2-category, and cannot use the expected result that
$3\Cat$ forms a weak 4-category.   Instead, we need to start
with a semistrict monoidal category $\C$, explicitly describe 
the objects, morphisms, and 2-morphisms of $\Z(\C)$, 
and then rather laboriously prove that it is indeed a semistrict
monoidal 2-category. 

In fact, it is natural to conjecture
a kind of `categorification' of the whole theory of quantum doubles.
For example, one should be able to start with
a `Hopf category' as defined by Crane and Frenkel \cite{CF} --- or, better,  
a `Hopf 2-algebra' --- and form the monoidal 2-category $\Reps(H)$ of its
representations on `2-vector spaces' \cite{KV,Yetter2}.  The monoidal
2-category $\Reps(H)$ should be equipped with a monoidal 2-functor to 
$2\Vect$ and satisfy various other conditions, and there 
should be a Tannaka--Krein theorem saying that, conversely,
such data determine a Hopf 2-algebra, unique up to equivalence.    
The center $\Z(\Reps(H))$ should thus
be a braided monoidal 2-category, and by Tannaka--Krein reconstruction
should determine a Hopf 2-algebra $DH$, the `quantum double'
of $H$.  Finally, one expects that this quantum double will be `quasitriangular'
in the sense defined by Crane and Frenkel \cite{CF}.  More ambitiously,
one might conjecture a similar correspondence between braided monoidal
$n$-categories and quasitriangular Hopf $n$-algebras for higher $n$.  We
shall not attempt to make these conjectures precise and prove them
here.  However, it is helpful to keep them in mind when considering
the applications of braided monoidal 2-categories to topology.

\subsection{Applications to 4-Dimensional TQFT}

Braided monoidal categories are especially interesting
because they give efficient procedures for constructing
tangle invariants and 3-dimensional
topological quantum field theories (TQFTs).  Braided monoidal 2-categories 
appear to have analogous applications to 2-tangle invariants and 
4-dimensional TQFTs.  As the TQFT applications are more intimately
related to the center construction, we begin with these.  To see
the patterns involved, it is helpful to consider first the rather 
trivial case of 2-dimensional TQFTs.

A 2-dimensional TQFT is a particular sort of symmetric monoidal functor $\F \maps 2\Cob \to \Vect$. Here the category $2\Cob$ has compact oriented 1-manifolds as objects and compact oriented cobordisms between them as morphisms, and it has a monoidal structure given by disjoint union. Similarly, the category $\Vect$ of finite-dimensional vector spaces and linear maps has a monoidal structure given by the usual tensor product.  In both cases these categories have a natural symmetric structure, as described in HDA and the references therein. The sphere with 3 open discs removed, or `trinion', can be 
thought of as a morphism in $2\Cob$:
\[            
    m \maps S^1 \cup S^1 \to S^1  , 
\]
and it gives rise to a product on the vector space $\F(S^1)$:
\[                \F(m) \maps \F(S^1)\otimes \F(S^1) \to \F(S^1)  .\]
One can easily check that this product is associative and 
commutative.  Similarly, the closed disc can be thought 
of as a morphism 
\[           i \maps \emptyset \to S^1, \]
which gives rise to a unit for the product on $\F(S^1)$:
\[                 \F(i) \maps {\mathbb{C}} \to \F(S^1)  .\]
Thus any 2-dimensional TQFT assigns to the circle a commutative
algebra.  

The true significance of this fact takes a bit of work to
unearth.  First, we can define a `commutative
monoid object' in any symmetric monoidal category to be an 
object $A$ equipped with a product and unit 
\[     m \maps A \otimes A \to A ,\qquad i \maps 1 \to A \]
satisfying analogs of the axioms for a commutative monoid.
In particular, $\F(S^1)$ is a commutative monoid object in $\Vect$,
that is, a commutative algebra.   However, this is really just
a corollary of the fact that
$S^1$ is a commutative monoid object in $2\Cob$, since
a symmetric monoidal functor takes commutative monoid objects to 
commutative monoid objects.   The real question is therefore, 
 \emph{why is $S^1$ a commutative monoid object in $2\Cob$?}

We shall not address this question directly.  Instead, note that whenever $A$ is a commutative monoid object in a symmetric monoidal category, $\hom(1,A)$ is a commutative monoid. Thus $\hom(\emptyset, S^1)$ is a commutative monoid.  Conversely,  understanding this commutative monoid should help us understand why $S^1$ is a commutative monoid object.  Moreover, by  following the patterns in Figure 1, we can learn something about the role of braided monoidal categories for  3-dimensional TQFTs, and braided monoidal 2-categories in 4-dimensional TQFTs.  

An element of $\hom(\emptyset, S^1)$ is an equivalence class of compact oriented 2-manifolds $M$ whose boundary has been identified with $S^1$. Alternatively, by fitting the circle inside a square in a standard way, we can think of $M$ as a 2-manifold with corners whose boundary is a square. Then, given $x, y \in \hom(\emptyset,S^1)$ we can define a `vertical' product $xy$ and a `horizontal' product $x \otimes y$ as shown in Figure 2.

\begin{center}
    \begin{figure}[h]
        \setlength{\unitlength}{0.01in}%
        \begin{picture}(260,80)(120,520)
        \thicklines
        \put(120,600){\line( 0,-1){ 80}}
        \put(120,520){\line( 1, 0){ 80}}
        \put(200,520){\line( 0, 1){ 80}}
        \put(200,600){\line(-1, 0){ 80}}
        \put(120,560){\line( 1, 0){ 80}}
        \put(300,520){\line( 1, 0){ 80}}
        \put(380,520){\line( 0, 1){ 80}}
        \put(380,600){\line(-1, 0){ 80}}
        \put(300,600){\line( 0,-1){ 80}}
        \put(340,600){\line( 0,-1){ 80}}
        \put(160,575){\makebox(0,0)[lb]{\raisebox{0pt}[0pt][0pt]{$ x$}}}
        \put(160,535){\makebox(0,0)[lb]{\raisebox{0pt}[0pt][0pt]{$ y$}}}
        \put(320,555){\makebox(0,0)[lb]{\raisebox{0pt}[0pt][0pt]{$ x$}}}
        \put(360,555){\makebox(0,0)[lb]{\raisebox{0pt}[0pt][0pt]{$ y$}}}
        \end{picture}
    \caption{Vertical and horizontal product in $\hom(\emptyset, S^1)$}
    \end{figure}
\end{center}

\noindent These products satisfy the exchange identity,
and taking $M$ to be the disc gives an element $1 \in \hom(\emptyset, S^1)$ 
that is a unit for both the horizontal and vertical product.
The Eckmann--Hilton argument then implies that
$\hom(\emptyset, S^1)$ is a commutative monoid.  We depict
this argument graphically in Figure 3.    

\begin{center}
\begin{figure}[h]
    \setlength{\unitlength}{0.009in}%
    \begin{picture}(640,80)(120,520)
    \thicklines
    \put(120,600){\line( 0,-1){ 80}}
    \put(120,520){\line( 1, 0){ 80}}
    \put(200,520){\line( 0, 1){ 80}}
    \put(200,600){\line(-1, 0){ 80}}
    \put(260,600){\line( 0,-1){ 80}}
    \put(260,520){\line( 1, 0){ 80}}
    \put(340,520){\line( 0, 1){ 80}}
    \put(340,600){\line(-1, 0){ 80}}
    \put(400,600){\line( 0,-1){ 80}}
    \put(400,520){\line( 1, 0){ 80}}
    \put(480,520){\line( 0, 1){ 80}}
    \put(480,600){\line(-1, 0){ 80}}
    \put(540,600){\line( 0,-1){ 80}}
    \put(540,520){\line( 1, 0){ 80}}
    \put(620,520){\line( 0, 1){ 80}}
    \put(620,600){\line(-1, 0){ 80}}
    \put(300,600){\line( 0,-1){ 80}}
    \put(260,560){\line( 1, 0){ 80}}
    \put(580,600){\line( 0,-1){ 80}}
    \put(540,560){\line( 1, 0){ 80}}
    \put(680,600){\line( 0,-1){ 80}}
    \put(680,520){\line( 1, 0){ 80}}
    \put(760,520){\line( 0, 1){ 80}}
    \put(760,600){\line(-1, 0){ 80}}
    \put(160,600){\line( 0,-1){ 80}}
    \put(400,560){\line( 1, 0){ 80}}
    \put(720,600){\line( 0,-1){ 80}}
    \put(140,555){\makebox(0,0)[lb]{\raisebox{0pt}[0pt][0pt]{$ x$}}}
    \put(175,555){\makebox(0,0)[lb]{\raisebox{0pt}[0pt][0pt]{$ y$}}}
    \put(280,535){\makebox(0,0)[lb]{\raisebox{0pt}[0pt][0pt]{$ x$}}}
    \put(315,575){\makebox(0,0)[lb]{\raisebox{0pt}[0pt][0pt]{$ y$}}}
    \put(435,535){\makebox(0,0)[lb]{\raisebox{0pt}[0pt][0pt]{$ x$}}}
    \put(435,575){\makebox(0,0)[lb]{\raisebox{0pt}[0pt][0pt]{$ y$}}}
    \put(555,575){\makebox(0,0)[lb]{\raisebox{0pt}[0pt][0pt]{$ y$}}}
    \put(595,535){\makebox(0,0)[lb]{\raisebox{0pt}[0pt][0pt]{$ x$}}}
    \put(735,555){\makebox(0,0)[lb]{\raisebox{0pt}[0pt][0pt]{$ x$}}}
    \put(695,555){\makebox(0,0)[lb]{\raisebox{0pt}[0pt][0pt]{$ y$}}}
    \put(275,575){\makebox(0,0)[lb]{\raisebox{0pt}[0pt][0pt]{ 1}}}
    \put(315,535){\makebox(0,0)[lb]{\raisebox{0pt}[0pt][0pt]{ 1}}}
    \put(595,575){\makebox(0,0)[lb]{\raisebox{0pt}[0pt][0pt]{ 1}}}
    \put(555,535){\makebox(0,0)[lb]{\raisebox{0pt}[0pt][0pt]{ 1}}}
    \put(225,555){\makebox(0,0)[lb]{\raisebox{0pt}[0pt][0pt]{$=$}}}
    \put(365,555){\makebox(0,0)[lb]{\raisebox{0pt}[0pt][0pt]{$=$}}}
    \put(505,555){\makebox(0,0)[lb]{\raisebox{0pt}[0pt][0pt]{$=$}}}
    \put(645,555){\makebox(0,0)[lb]{\raisebox{0pt}[0pt][0pt]{$=$}}}
    \end{picture}
    
    \caption{The Eckmann--Hilton argument}
    \end{figure}
\end{center}

The appearance of the Eckmann--Hilton argument here suggests that
we really have a 2-category with one object and one 1-morphism
on our hands.  Now, the `extended TQFT hypothesis' in HDA 
suggests that the best way to understand $n$-dimensional TQFTs
is in terms of a weak $n$-category $\C_{n,\infty}$ whose objects are
0-manifolds, whose morphisms are equivalence classes of 1-manifolds 
with boundary, whose 2-morphisms are equivalence classes of
2-manifolds with corners, and so on, each $(j+1)$-morphism
being a kind of cobordism between $j$-morphisms.  
(Of course these manifolds should be compact and oriented;
in general they should also be `framed', but
here we neglect this subtlety.)    Making this hypothesis precise
would require a general definition of weak $n$-categories, and
also some careful differential topology.   Even in its current
vague form, though, it sheds some light on the situation at hand.  
$\C_{n,\infty}$ 
should have a distinguished object $\ast$, the positively oriented point.  
The 1-morphism $1_\ast$ should then correspond to the closed unit interval.
When $n = 2$, $\hom(1_\ast,1_\ast)$ should then be the set of all cobordisms
from the interval to itself.  These are just equivalence
classes of 2-manifolds with corners whose boundary is the square!  
Thus $\hom(1_\ast,1_\ast)$ is isomorphic to
$\hom(\emptyset, S^1)$, but now the commutative monoid
structure has a purely algebraic explanation: there is a 2-category having 
one object $\ast$, one 1-morphism $1_\ast$, and the set $\hom(1_\ast,1_\ast)$
as its 2-morphisms.   Understand the isomorphism between
$\hom(\emptyset,S^1)$ and $\hom(1_\ast,1_\ast)$ in purely
algebraic terms remains an interesting challenge; the solution 
will probably involve
the theory of duality in $n$-categories.  

Similarly, in the study of 3-dimensional
TQFTs we expect to have a 3-category $\C_{3,\infty}$, and sitting
inside this there should be a 3-category with one object $\ast$, one
1-morphism $1_\ast$, and the category $\hom(1_\ast, 1_\ast)$ as
its 2-morphisms and 3-morphisms.   This category should thus 
be a braided monoidal category whose objects 
are 2-manifolds with corners having a square
as boundary, and whose morphisms are cobordisms between
these.  Likewise, in the 4-dimensional case $\hom(1_\ast, 1_\ast)$ would
be a braided monoidal 2-category, and so on.  

In fact, results along these lines already appear  in the literature
in the cases of dimensions 3 and 4, but in terms of  
$\hom(\emptyset,S^1)$ rather than $\hom(1_\ast,1_\ast)$.
This is less natural algebraically,
but simpler topologically, because the theory of cobordisms between
manifolds with corners is not well developed.
So far, the clearest description of $\hom(\emptyset,S^1)$
as a braided monoidal category in dimension 3 and a braided
monoidal 2-category in dimension 4 has been given by Crane and Yetter \cite{CY}.
There are many interesting projects left to do, however.
For example, in dimension 3 it should be possible to use existing results of Kerler 
\cite{Kerler} and others to obtain a presentation of $\hom(\emptyset, S^1)$ 
as a braided monoidal category, and to compare the answer
to what one would predict using the extended TQFT hypothesis.
This presentation should explain the already known conditions required
to construct 3-dimensional TQFTs, such as Chern--Simons theory,
which associate a braided monoidal category to the circle
\cite{Crane,RT2}.  
In dimension 4 one still needs to carefully check whether $\hom(\emptyset,
S^1)$ meets our definition of a braided monoidal 2-category, and then if
possible obtain a presentation of it.   This may allow the construction
of 4-dimensional TQFTs from braided monoidal 2-categories meeting certain
conditions.  If so, our center construction may serve as a source of
4-dimensional TQFTs.  

\subsection{Applications to 2-Tangles}
\label{2tangles}

Tangles can be regarded as certain equivalence classes of 1-manifolds with
boundary embedded in  $[0,1]^3$, possibly equipped with extra structure
such as an orientation or framing.  Tangles are important because they make
clear the relation between knot theory and braided monoidal categories.  
For example, framed oriented tangles form 
the `free balanced braided monoidal category on one object' 
\cite{FY,JS,Turaev,Yetter}, and this fact permits the construction 
of knot invariants from the categories of representations of 
quantum groups and other quasitriangular Hopf algebras \cite{RT}.

The `tangle hypothesis' of HDA suggests that this is part of a more general 
relationship between `$k$-tangles in $(n+k)$ dimensions' and 
$k$-tuply monoidal $n$-categories.  A $k$-tangle in $(n+k)$
dimensions is something like an isotopy equivalence class of 
$k$-manifolds with corners embedded in $[0,1]^{n+k}$.   The tangle hypothesis
proposes that these may be described algebraically using a specific $k$-tuply
monoidal $n$-category $\C_{n,k}$, which has the cobordism $n$-category
$\C_{n,\infty}$ as a limiting case.  

A very interesting example is the case of 2-tangles in 4 dimensions: 
$n = 2$, $k = 2$.    Topologists have already studied these 2-tangles,
and the work of Carter and Saito \cite{CS} strongly suggests that
they form a braided monoidal 2-category.  In fact, Fischer \cite{Fischer}
claims to have already shown this.  His work is unfortunately rather
unclear, but Kharlamov and Turaev \cite{KT} have begun to redo it more carefully.
It should also be re-evaluated in the light of our definition of braided 
monoidal 2-category.  One would eventually like to construct 2-tangle
invariants from certain braided monoidal 2-categories, such as the
category of representations of quasitriangular Hopf 2-algebras.
Our center construction is a small step in this direction.

We cannot conclude this introduction without a word or two
about the Zamolodchikov tetrahedron equation.  In a braided monoidal
category, the braiding automatically satisfies the Yang--Baxter
equation.  In other words, given objects $A,B,C$, the following
diagram commutes:

\[
\begin{tikzcd}[column sep = large]
    &
    B \otimes A \otimes C
    \arrow[r, "B \otimes R_{A,C}"]
    &
    B \otimes C \otimes A
    \arrow[dr, "R_{B,C} \otimes A"]
    \\
    A \otimes B \otimes C
    \arrow[ur, "R_{A,B} \otimes C"]
    \arrow[dr, "A \otimes R_{B,C}", swap]
    &&&
    C \otimes B \otimes A
    \\&
    A \otimes C \otimes B
    \arrow[r, "R_{A,C} \otimes B", swap]
    &
    C \otimes A \otimes B
    \arrow[ur, "C \otimes R_{A,B}", swap]
\end{tikzcd}\]

\noindent  In the theory of tangles this corresponds to the following
equation between tangles:

\begin{figure}[h]
\[
\begin{tikzpicture}
    \node[coordinate] (a0) at (0,0) {};
    \node[coordinate] (b0) at (1,0) {};
    \node[coordinate] (c0) at (2,0) {};
    \node (ab1) at (0.5,0.7) {};
    \node (c1)  at (2,0.7)   {};
    \node (c2)  at (2,1.4)   {};
    \node (bc3) at (1.5,2.1) {};
    \node (c3)  at (2,2.1)   {};
    \node (ab5) at (0.5,3.5) {};
    \node[coordinate] (a6) at (0,4.2) {};
    \node[coordinate] (b6) at (1,4.2) {};
    \node[coordinate] (c6) at (2,4.2) {};
    
    \draw (a0) .. controls (ab1) and (c3) .. (c6);
    \draw (b0) to (ab1) to [bend left] (ab5) to (b6);
    \draw (c0) .. controls (c1) and (c2) .. (bc3)
          (bc3) to (ab5) (ab5) to (a6);
\end{tikzpicture}
\quad \raisebox{2.1cm}{=} \quad
\begin{tikzpicture}
    \node (bc1) at (1.5,0.7) {};
    \node (a3)  at (0,2.1)   {};
    \node (ab3) at (0.5,2.1) {};
    \node (a4)  at (0,2.8)   {};
    \node (a5)  at (0,3.5)   {};
    \node (bc5) at (1.5,3.5) {};
    
    \draw (c6) .. controls (bc5) and (a3) .. (a0);
    \draw (b6) to (bc5) to [bend left] (bc1) to (b0);
    \draw (a6) .. controls (a5) and (a4) .. (ab3)
          (ab3) to (bc1) (bc1) to (c0);
\end{tikzpicture}
\]
\caption{The Yang--Baxter equation}
\end{figure}
\vskip1em

In a braided monoidal 2-category, the Yang--Baxter equation
holds only up to a 2-isomorphism.  Topologically, this 2-isomorphism
corresponds to a 2-tangle which intersected with $\{0\} \times [0,1]^3 \subset
[0,1]^4$ looks like the left side of Figure 4, and which intersected
with $\{1\} \times [0,1]^3$ looks like the right side of Figure 4.  

In Kapranov and Voevodsky's theory \cite{KV} there are in fact
two distinct such 2-isomorphisms, $S^{\pm}_{A,B,C}$, corresponding
to two distinct \emph{proofs} of the Yang--Baxter equation in a
braided monoidal category.  However, they give the same 
2-tangle.  There is also a deep relationship between $n$-category theory
and homotopy theory, described in HDA and the references therein, and using
this, Breen \cite{Breen} has deduced that the condition $S^+ = S^-$ 
should hold.

These facts constitute topological evidence that in the correct
definition of a braided monoidal category, there should be an extra
coherence law asserting that $S^+ = S^-$.  We also find algebraic
evidence for this, as follows.  It follows heuristically from our rough
definition of center that a $k$-tuply monoidal $n$-category should embed
canonically in its center when $\C$ happens to be already $(k+1)$-tuply
monoidal.  More precisely, if $\C$ is a $(k+1)$-tuply monoidal
$n$-category and $\C_0$ is the underlying $k$-tuply monoidal
$n$-category, there should be a faithful $(k+1)$-tuply monoidal
$n$-functor from $\C$ to $\Z(\C_0)$.  For example, the center of a set
$S$ works out to be the monoid ${\rm End}(S)$, but when $S$ happens
already to be a monoid, there is a natural embedding $S \hookrightarrow
{\rm End}(S)$ given by the left action of $S$ on itself.  Similarly, a
monoid equals its center when it is commutative, and a monoidal category
naturally embeds in its center when it is braided \cite{JS,Kassel}.  The
third main goal of this paper is to show that a monoidal 2-category $\C$
embeds into $\Z(\C)$ when $\C$ happens to be braided.  However, for any
monoidal 2-category $\C$ it turns out that $S^+ = S^-$ in $\Z(\C)$.
Thus we can only achieve our goal if our definition of braided monoidal
2-category includes a coherence law saying that $S^+ = S^-$.  (It is
worth noting that all our results except Theorem~\ref{embed} hold
without this extra coherence law.)

Finally, if $S = S^+ = S^-$, Kapranov and Voevodky's work \cite{KV}
implies that the 2-morphisms $S_{A,B,C}$ satisfy an equation of their
own, the Zamolodchikov tetrahedron equation.  This is the
higher-dimensional analogue of the Yang--Baxter equation, and it plays an
important role in the theory of 2-tangles.  Pictures of the
Zamolodchikov tetrahedron equation in terms of 2-tangles can be found in
the work of Carter and Saito \cite{CS}.  Kapranov and Voevodsky, who do
not assume $S^+ = S^-$, write down 8 different versions of the
Zamolodchikov equation and claim that these all follow from their
definition of a braided monoidal 2-category.  In our framework there is
only one Zamolodchikov equation.

%% file: 2definitions.tex
\section{Definitions}
We begin by defining semistrict monoidal 2-categories and semistrict braided 
monoidal 2-categories.    Following traditional practice among category
theorists \cite{Gray,KS}, we use `2-category' to mean what Kapranov and
Voevodsky \cite{KV} call a strict 2-category, and `2-functor' to mean
what Kapranov and Voevodsky call a strict 2-functor.  Composition of
1-morphisms, the horizontal composition of a 1-morphism and a 2-morphism
(in either order) and the horizontal composition of
2-morphisms is denoted by $\circ$ or simply juxtaposition.  
Vertical composition of 2-morphisms is denoted by $\cdot$.
We use the ordering in which, for example, the composite of 
$f \maps A \to B$ and $g \maps B \to C$ 
is denoted $f\circ g$.  

We use $\C \gtimes \D$ to denote Gordon, Power, and Street's
\cite{GPS} `Gray' tensor product of the 2-categories 
$\C$ and $\D$.  This differs
from Gray's original version \cite{Gray} in being the `pseudo' rather than 
the `lax' weakening of the Cartesian product.    For readers
unfamiliar with these distinctions, let us simply recall that  
given a 1-morphism $f \maps A \to A'$ in $\C$ and
a 1-morphism $g \maps B \to B'$ in $\D$, the Cartesian
product $\C \times \D$ contains a commuting square
 
\[
\begin{tikzcd}
    {(A,B)}
    \arrow[r, "f \times 1"]
    \arrow[d, "1 \times g", swap]
    &
    {(A,B)}
    \arrow[d, "1 \times g"]
    \\
    {(A,B)}
    \arrow[r, "f \times 1", swap]
    &
    {(A',B')}
\end{tikzcd}\]

Following the `lax' approach to weakening, which consists of replacing
equations by morphisms,
Gray's original product of $\C$ and $\D$ instead contains a square commuting
only up to a specified 2-morphism:

\[
\begin{tikzcd}
    {(A,B)}
    \arrow[r, "f \gtimes 1"]
    \arrow[d, "1 \gtimes g", swap]
    &
    {(A,B)}
    \arrow[dl, phantom, "\Downarrow \gamma_{f,g}"]
    \arrow[d, "1 \gtimes g"]
    \\
    {(A,B)}
    \arrow[r, "f \gtimes 1", swap]
    &
    {(A',B')}
\end{tikzcd}\]

Following the `pseudo' approach, which consists of replacing equations
by isomorphisms (or, more generally, equivalences),
Gordon, Power, and Street additionally require 
$\gamma_{f,g}$ to be an isomorphism.   We use their version 
of the Gray tensor product as part of a systematic
adherence to the `pseudo' approach.

The category $2\Cat$ with
2-categories as objects and 2-functors as morphisms
becomes a monoidal category $(2\Cat, \gtimes, \I)$
when equipped with the Gray 
tensor product and the unit object $\I$, the 2-category
with one object, one morphism
and one 2-morphism.   This monoidal category is symmetric, with the symmetry
\[               S_{\C,\D} \maps \C \gtimes \D \to \D \gtimes \C \]
given by:
\[ 
\begin{array}{ccc}
(B,A)  \mapsto (B,A) & 
(f,1) \mapsto (1,f) &
(1,g) \mapsto (g,1)  \\
(\alpha,1) \mapsto (1,\alpha) & 
(1,\beta) \mapsto (\beta,1) & 
\gamma_{f,g} \mapsto \gamma_{g,f}^{-1}
\end{array} \]

\subsection{Semistrict Monoidal 2-Categories}

Since $2\Cat$ is monoidal when equipped with the Gray tensor product, we may use enriched category theory \cite{Kelly} to efficiently define semistrict 3-categories and monoidal 2-categories:

\begin{defn}
    A \define{semistrict 3-category} is a category enriched over $(2\Cat,\gtimes,\I)$.
\end{defn}

\begin{defn}
    A \define{semistrict monoidal 2-category} is a semistrict 3-category with one object.  
\end{defn}
 
Gordon, Power and Street \cite{GPS} have given a definition of   `weak' 3-categories, or `tricategories', seemingly more general than that of semistrict 3-categories, and indeed intended to be `maximally general' in some sense.  For example, associativity and identity laws hold as equations in a semistrict 3-category, but only hold up to  specified equivalence in a weak one.  However, these authors have shown that every weak 3-category is equivalent in a precise sense (`triequivalence') to a semistrict one, so for many purposes semistrict 3-categories are `sufficiently general'. Defining a weak monoidal 2-category to be a weak 3-category with one object, it follows  from their proof that any one of these is triequivalent to a semistrict monoidal 2-category.  So again, while not maximally general, semistrict
monoidal 2-categories are sufficiently general for many purposes.

Often we shall think of a semistrict
monoidal 2-category as a 2-category with extra structure.  
More precisely, if $\tilde \C$ is a semistrict 3-category with one object $\ast$, let 
$\C = \hom(\ast, \ast)$.  This is a 2-category equipped with
a 2-functor 
\[                \otimes: \C \gtimes \C \to \C  \]
coming from composition in $\tilde \C$, as well as a functor
$i \maps \I \to \C$ coming from the identity of $\ast$ in $\tilde \C$.

\begin{lma} \label{lem1} Suppose $\tilde \C$ is a semistrict 3-category with
one object, and let $(\C, \otimes, i)$ be defined as above.
Then $\C$ is a 2-category, $\otimes \maps \C \gtimes \C \to \C$ 
and $i \maps \I \to \C$ are 2-functors, and the following diagrams
commute:
\begin{enumerate}
\item Associativity:

\[
\begin{tikzcd}
    \C \gtimes \C \gtimes \C
    \arrow[r, "\otimes \gtimes \C"]
    \arrow[d, "\C \gtimes \otimes", swap]
    &
    \C \gtimes \C
    \arrow[d, "\otimes"]
    \\
    \C \gtimes \C
    \arrow[r, "\otimes", swap]
    &
    \C
\end{tikzcd}\]

\item Unit law:

\[
\begin{tikzcd}[column sep = small]
    \I \gtimes \C
    \arrow[rr, "i \gtimes \C"]
    \arrow[dr, "\cong", swap]
    &&
    \C \gtimes \C
    \arrow[dl, "\otimes"]
    \\&
    \C
\end{tikzcd}
\quad\quad\quad
\begin{tikzcd}[column sep = small]
    \C \gtimes \I
    \arrow[rr, "\C \gtimes i"]
    \arrow[dr, "\cong", swap]
    &&
    \C \gtimes \C
    \arrow[dl, "\otimes"]
    \\&
    \C
\end{tikzcd}\]

\end{enumerate}

Conversely, for any $(\C,\otimes,i)$ with these properties, 
there is a unique semistrict 3-category
$\tilde \C$ with one object from which $(\C,\otimes,i)$ arises
as above.  \end{lma}

\begin{proof}
    This is a straightforward consequence of the definition of semistrict 3-categories as categories enriched over $2\Cat$ with its Gray tensor product.
\end{proof}

There is thus no harm in thinking of a semistrict monoidal 2-category as a triple $(\C,\otimes,i)$ satisfying the associativity and unit law conditions of Lemma~\ref{lem1}.   Since the 2-functor $i$ is determined by the object of  $\C$ obtained by applying $i \maps \I \to \C$ to the one object in $\I$,  we can also think of a semistrict monoidal 2-category as a triple $(\C,\otimes,I)$.

One may further unpack our definition of a 
semistrict monoidal 2-category and obtain the same explicit list of
operations and 
laws that Kapranov and Voevodsky take as their definition \cite{KV}.  Here
the standard machinery of 2-categorical commutative diagrams becomes
very handy \cite{KS}.  In what follows we write $\otimes_{f,g}$ for the 
2-morphism $\otimes(\gamma_{f,g})$ in $\C$.   

\begin{lma}   
\label{mon}
A semistrict monoidal 2-category consists of a 2-category 
$\C$ together with:
\begin{enumerate}
\item An object $I \in \C$.
\item For any two objects $A, B$ in $\C$, an object $A \otimes B$ in $\C$.

\item For any 1-morphism $f\maps A \to A'$ and any object $B \in \C$ a 
1-morphism $f \otimes B \maps A \otimes B \to A' \otimes B$.

\item For any 1-morphism $g\maps B \to B'$ and any object $A \in \C$ a 1-morphism
$A \otimes g \maps A \otimes B \to A \otimes B'$.

\item For any object $B \in \C$ and any 2-morphism 
$\alpha \maps f \Rightarrow f'$ a 2-morphism 
$\alpha \otimes B\maps f\otimes B \Rightarrow f' \otimes B$.

\item For any object $A \in \C$ and any 2-morphism 
$\beta \maps g \Rightarrow g' $ a 2-morphism
$A \otimes \beta \maps A \otimes g \Rightarrow A \otimes g'$.

\item For any two 1-morphisms $f\maps A \to A'$ and $g\maps B \to B'$ a 2-isomorphism

\[
\begin{tikzcd}[row sep = huge]
    A \otimes B
    \arrow[r, "A \otimes g"]
    \arrow[dr, phantom, "\Downarrow \otimes_{f,g}"]
    \arrow[d, "f \otimes B", swap]
    &
    A \otimes B'
    \arrow[d, "f \otimes B'"]
    \\
    A' \otimes B
    \arrow[r, "A' \otimes g", swap]
    &
    A' \otimes B'
\end{tikzcd}\]

\end{enumerate}
Moreover, these data must satisfy the following conditions.
\begin{itemize}
\item[(i)] For any object $A \in \C$ we have $A \otimes \stri \; \maps\C \to \C$ and
$ \stri \otimes A \maps \C \to \C$ are 2-functors.
\item[(ii)]  For $x$ any object, morphism or 2-morphism of $\C$ we have
$x \otimes I = I \otimes x = x $.
\item[(iii)]  For $x$ any object, morphism or 2-morphism of $\C$, and 
for all objects $A,B \in \C$ we
have $A \otimes (B \otimes x) = (A \otimes B) \otimes x$,
$A \otimes (x \otimes B) = (A \otimes x) \otimes B$ and 
$x \otimes (A \otimes B) = (x \otimes A) \otimes B$.  
\item[(iv)]
For any 1-morphisms $f\maps A \to A'$, $g\maps B \to B'$ and 
$h \maps C \to C'$ in $\C$ 
we have $\bigotimes_{A \otimes g,h} = A \bigotimes \otimes_{g,h}$, 
$\bigotimes_{f y\otimes B,h} = \bigotimes_{f, B \otimes h}$ and 
$\bigotimes_{f,g \otimes C} = \bigotimes_{f,g} \otimes C$.
\item[(v)]  For any objects $A,B \in \C$ we have $1_A \otimes B = 
A \otimes 1_B = 1_{A \otimes B}$, and for any 1-morphisms $f \maps 
A \to A'$, $g \maps B \to B'$ in $\C$ we have $\bigotimes_{1_A,g} = 
1_{A \otimes g}$  and $\bigotimes_{f,1_B} = 1_{f \otimes B}$.  
\item[(vi)] For any 1-morphism $f: A \to A'$, any 1-morphisms
$g,g' \maps B \to B'$, and any 2-morphism
$\beta \maps g \Rightarrow g' $ the following diagram commutes: 
\smallskip

\[
\begin{tikzcd}[row sep = large]
    A \otimes B
    \arrow[ddr, phantom, "\Downarrow \otimes_{f,g}", pos = 0.35]
    \arrow[rr, bend left]
    \arrow[rr, bend right]
    \arrow[dd]
    &
    \zwei A \otimes \beta
    &
    A \otimes B'
    \arrow[dd]
    \\\\
    A' \otimes B
    \arrow[rr, dashed, bend left]
    \arrow[rr, bend right]
    &
    \zwei A' \otimes \beta
    &
    A' \otimes B'
    \arrow[uul, phantom, "\Downarrow \otimes_{f,g'}", pos = 0.35]
\end{tikzcd}\]

\smallskip

\item[(vii)]  For any 1-morphism $g: B \to B'$, any 1-morphisms $f,f'
\maps A \to A'$, and any 2-morphism
$\alpha \maps f \Rightarrow f'$, the following diagram commutes: \smallskip

\[
\begin{tikzcd}[row sep = large]
    A \otimes B
    \arrow[ddr, phantom, "\Uparrow \otimes_{f,g}", pos = 0.35]
    \arrow[rr, bend left]
    \arrow[rr, bend right]
    \arrow[dd]
    &
    \zwei \alpha \otimes B
    &
    A' \otimes B
    \arrow[dd]
    \\\\
    A \otimes B'
    \arrow[rr, dashed, bend left]
    \arrow[rr, bend right]
    &
    \zwei \alpha \otimes B'
    &
    A' \otimes B'
    \arrow[uul, phantom, "\Uparrow \otimes_{f',g}", pos = 0.35]
\end{tikzcd}\]

\smallskip

\item[(viii)]
 For any 1-morphisms $f\maps A \to A'$, $g\maps B \to B'$ and $g'\maps
B' \to B''$  the 2-isomorphism $\bigotimes_{f,gg'}$ coincides with the
pasting of  $\bigotimes_{f,g}$ and $\bigotimes_{f,g'}$ as in the
following diagram. \medskip 

\[
\begin{tikzcd}[column sep = small]
    A \otimes B
    \arrow[rr]
    \arrow[dd]
    &&
    A \otimes B'
    \arrow[rr]
    \arrow[dd]
    \arrow[ddll, phantom, "\zwei \otimes_{f,g}"]
    \arrow[ddrr, phantom, "\zwei \otimes_{f,g'}"]
    &&
    A \otimes B''
    \arrow[dd]
    \\\\
    A' \otimes B
    \arrow[rr]
    &&
    A' \otimes B'
    \arrow[rr]
    &&
    A' \otimes B''
\end{tikzcd}
\]                                      
\medskip         
\noindent
 For any 1-morphisms $f\maps A \to A'$, $f'\maps A' \to A''$ and $g\maps B \to B'$
 the 2-isomorphism $\bigotimes_{ff',g}$ coincides with the pasting of 
$\bigotimes_{f,g}$ and $\bigotimes_{f,g'}$  in a similar way.

\end{itemize}
\end{lma}

\begin{proof} This is a straightforward verification.  In particular,
conditions (v), (vi) and (vii) come from the coherence laws satisfied by
$\gamma_{f,g}$ in the Gray tensor product. \end{proof}
Note that condition $(viii)$ and the invertibility of the $2$-morphism 
$\otimes_{f,g}$ imply that $\otimes_{1_A,g} = 1_g$ and 
$\otimes_{f,1_B} = 1_f$, for any $f : A \to A'$ and any $g : B \to B'$.  

\subsection{Semistrict Braided Monoidal 2-Categories}

To efficiently define braided monoidal 2-categories it is
useful to exploit the fact that $(2\Cat,\gtimes,\I)$ is closed, i.e.,
enriched over itself \cite{GPS}.  Put more explicitly,
what this means is that $2\Cat$ 
can be regarded as a semistrict 3-category having small 
2-categories as objects, 2-functors as morphisms,
`pseudonatural transformations' as 2-morphisms, and 
`modifications' as 3-morphisms \cite{BG,KS}.  A pseudonatural
transformation $T$ 
between 2-functors $\F,\G \maps \C \to \D$ assigns to each 
object $A \in \C$ a morphism $T_A \maps \F(A) \to \G(A)$
which satisfies the definition of a natural transformation only {\it up
to a specified isomorphism}.  
Thus, $T$ also assigns to each morphism $f \maps A \to B$ in
$\C$ a 2-isomorphism $T_f$ as follows: 

\[
\begin{tikzcd}[row sep = huge, column sep = large]
    \F(A)
    \arrow[r, "\F(f)"]
    \arrow[d, swap, "T_A"]
    \arrow[dr, phantom, "\zwei T_f"]
    &
    \F(B)
    \arrow[d, "T_B"]
    \\
    \G(A)
    \arrow[r, swap, "\G(f)"]
    &
    \G(B)
\end{tikzcd}
\]

\noindent  These 2-morphisms $T_f$ must in turn satisfy some
equational laws of their own.  First, for any identity morphism
$1_A \maps A \to A$, we require $T_{1_A} = 1_{T_A}$.  Second,
given a composable pair of morphisms $f \maps A \to B$, $g \maps B \to C$, 
the 2-morphism $T_{fg}$ is given by the following pasting:
\medskip

\[
\begin{tikzcd}[row sep = huge, column sep = large]
    \F(A)
    \arrow[r]
    \arrow[d]
    \arrow[dr, phantom, "\zwei T_f"]
    &
    \F(B)
    \arrow[r]
    \arrow[d]
    \arrow[dr, phantom, "\zwei T_g"]
    &
    \F(C)
    \arrow[d]
    \\
    \G(A)
    \arrow[r]
    &
    \G(B)
    \arrow[r]
    &
    \G(C)
\end{tikzcd}
\]
\bigskip    
  
\noindent Third, given morphisms $f,f' \maps A \to B$ and a
2-morphism $\alpha \maps f \Rightarrow f'$, the following diagram commutes:

\[
\begin{tikzcd}[row sep = large]
    \F(A)
    \arrow[ddr, phantom, "\zwei T_{f'}", pos = 0.25]
    \arrow[rr, bend left]
    \arrow[rr, bend right]
    \arrow[dd]
    &
    \F(\alpha)
    &
    \F(B)
    \arrow[dd]
    \\\\
    \G(A)
    \arrow[rr, dashed, bend left]
    \arrow[rr, bend right]
    &
    \G(\alpha)
    &
    \G(B)
    \arrow[uul, phantom, "\zwei T_{f}", pos = 0.25]
\end{tikzcd}\]

\noindent Given two pseudonatural transformations $S,T \maps \F \Rightarrow \G$, 
a modification $\alpha$ from $S$ to $T$ assigns to each object $A \in \C$ a 2-morphism $\alpha_A \maps S_A \Rightarrow T_A$.
Moreover, for any morphism $F \maps A \to B$, the following
diagram is required to commute:

\[
\begin{tikzcd}[row sep = large]
    \F(A)
    \arrow[ddr, phantom, "\Uparrow T_{f}", pos = 0.25]
    \arrow[rr, bend left]
    \arrow[rr, bend right]
    \arrow[dd]
    &
    \zwei \alpha_A
    &
    \G(A)
    \arrow[dd]
    \\\\
    \F(B)
    \arrow[rr, dashed, bend left]
    \arrow[rr, bend right]
    &
    \zwei \alpha_B
    &
    \G(B)
    \arrow[uul, phantom, "\Uparrow S_{f}", pos = 0.25]
\end{tikzcd}\]

As explained in the introduction, in an $n$-category
the notion of `isomorphism' can 
be weakened to a recursively defined notion of `equivalence'.
In the case of $2\Cat$ this gives the following concepts.
A modification $\alpha$ from the pseudonatural transformation
$S$ to the pseudonatural transformation $T$ is `invertible' if 
there is a modification $\alpha^{-1}$ from $T$ to $S$ such
that $\alpha\alpha^{-1} = 1_S$ and $\alpha^{-1}\alpha = 1_T$.
A pseudonatural transformation $T$ from $\F$ to $\G$ is a
`pseudonatural equivalence' if there is a pseudonatural transformation 
$\overline T\maps \G \to F$ and invertible modifications
\[     \alpha_1 \maps T\overline{ T} \to 1_\F,\qquad \alpha_2 \maps
\overline{ T} T \to 1_\G .\]
There is a similar notion at the level of 2-functors, but we will
not need it.

Every semistrict monoidal 2-category has a second, `opposite' tensor product:
\begin{lma}
Suppose $(\C,\otimes,I)$ is a semistrict monoidal 2-category.  Then 
$(\C,\otimes^\op,I)$ is also a semistrict monoidal 2-category, where 
$ \otimes^\op = S_{\C,\C} \circ \otimes $.  
\end{lma}
\begin{proof} Straightforward.  \end{proof}

There is an analogous opposite tensor product for
strict monoidal categories, and a strict braided monoidal category
is just a strict monoidal category equipped with a natural isomorphism
$R\maps \otimes \Rightarrow \otimes^\op$, the `braiding', 
such that the following triangles commute:

\[
\begin{tikzcd}[row sep = huge, column sep = tiny]
    A \otimes X \otimes Y
    \arrow[rr, "R_{A,X \otimes Y}"]
    \arrow[dr, swap, "R_{A,X} \otimes Y"]
    &&
    X \otimes Y \otimes A
    \\
    &
    X \otimes A \otimes Y
    \arrow[ur, swap, "X \otimes R_{A,Y}"]
\end{tikzcd}
\]
\bigskip

\[
\begin{tikzcd}[row sep = huge, column sep = tiny]
    X \otimes Y \otimes A
    \arrow[rr, "R_{X \otimes Y,A}"]
    \arrow[dr, swap, "X \otimes R_{Y,A}"]
    &&
    A \otimes X \otimes Y
    \\
    &
    X \otimes A \otimes Y
    \arrow[ur, swap, "R_{X,A} \otimes Y"]
\end{tikzcd}
\]

The definition of a semistrict braided monoidal 2-category is very similar.
However, instead of a strict monoidal category, one starts with a semistrict
monoidal 2-category.  Instead of the braiding being a natural
transformation, it is a pseudonatural equivalence.   Instead of
the equations above holding `on the nose', they hold up
to specified invertible modifications.  Finally, these modifications
must satisfy 3 new coherence laws discovered by Kapranov and Voevodsky,
together with the equation $S^+ = S^-$ discussed in Section~\ref{2tangles}.

In all that follows, in diagrams we sometimes denote the tensor product of
objects simply by juxtaposition.  We also label some clauses in the definition 
using the `hieroglyphic' notation invented by Kapranov and Voevodsky.

\begin{defn} \label{bm2} 
A \define{braided monoidal 2-category} $(\C, \otimes, I, R, \tilR_{(\stri|\stri, \stri)}, \tilR_{(\stri, \stri | \stri)})$ consists of:
\begin{enumerate}
\item
A semistrict monoidal 2-category $(\C,\otimes,1)$
\item
A pseudonatural equivalence 
$ R \maps \otimes \Rightarrow \otimes^\op$ 
\item
Two invertible modifications $\tilR_{(\stri|\stri ,\stri)}$ and 
$\tilR_{(\stri,\stri | \stri)}$, giving for any objects $A,B,C \in \C$
the 2-isomorphisms    

\[
\begin{tikzcd}[row sep = huge, column sep = tiny]
    A \otimes B \otimes C
    \arrow[rr, "R_{A,B \otimes C}"]
    \arrow[dr, swap, "R_{A,B} \otimes C"]
    & \phantom{0} &
    B \otimes C \otimes A
    \\
    &
    B \otimes A \otimes C
    \arrow[ur, swap, "B \otimes R_{A,C}"]
    \arrow[u, phantom, "\Uparrow \tilR_{(A|B,C)}", pos=0.7]
\end{tikzcd}
\qquad
\begin{tikzcd}[row sep = huge, column sep = tiny]
    A \otimes B \otimes C
    \arrow[rr, "R_{A \otimes B,C}"]
    \arrow[dr, swap, "A \otimes R_{B,C}"]
    & \phantom{0} &
    C \otimes A \otimes B
    \\
    &
    A \otimes C \otimes B
    \arrow[ur, swap, "R_{A,C} \otimes B"]
    \arrow[u, phantom, "\Uparrow \tilR_{(A,B|C)}", pos=0.7]
\end{tikzcd}
\]

\end{enumerate} \medskip

These data must satisfy the following conditions.  First,    
for all objects $A,B,C,D \in \C$ the following diagrams commute:

\vbox{
$((\bullet\otimes \bullet \otimes \bullet )\otimes \bullet)$

\[
\begin{tikzcd}
    &
    DABC
    \arrow[dddl, phantom, "1.", pos = 0.15, xshift=0.5em]
    \arrow[dddr, phantom, "2.", pos = 0.15, xshift=-0.25em]
    & \phantom{4}
    \\\\
    ABCD
    \arrow[uur]
    \arrow[dr]
    \arrow[rr, dashed]
    \arrow[uurr, phantom, "4.", pos = 0.06]
    &&
    ADBC
    \arrow[uul]
    \arrow[dll, phantom, "3.", pos = 0.1]
    \\
    \phantom{1} &
    ABDC
    \arrow[ur]
    \arrow[uuu]
    & \phantom{2}
\end{tikzcd}
\]

\[ \begin{array}{ll}
1. \:  \tilR_{(A \otimes B,C|D)} & 
2. \: =\:  \tilR_{(A,B|D)} \otimes C \\
3. \: = \:  A \otimes \tilR_{(B,C|D)} & 
4.\: =\:  \tilR_{(A,B\otimes C|D)}
\end{array}  \]
}

\vbox{
$(\bullet\otimes (\bullet \otimes \bullet\otimes \bullet))$ \medskip

\[
\begin{tikzcd}
    &
    BCDA
    \arrow[dddl, phantom, "1.", pos = 0.15, xshift=0.5em]
    \arrow[dddr, phantom, "2.", pos = 0.15, xshift=-0.25em]
    & \phantom{4}
    \\\\
    ABCD
    \arrow[uur]
    \arrow[dr]
    \arrow[rr, dashed]
    \arrow[uurr, phantom, "4.", pos = 0.06]
    &&
    BCAD
    \arrow[uul]
    \arrow[dll, phantom, "3.", pos = 0.1]
    \\
    \phantom{1} &
    BACD
    \arrow[ur]
    \arrow[uuu]
    & \phantom{2}
\end{tikzcd}
\]

\[ \begin{array}{ll}
1.\: =\:  \tilR_{(A| B,C\otimes D)} & 
2.\: =\:  B \otimes \tilR_{(A|C,D)}  \\
3.\: =\:  \tilR_{(A|B,C)}  \otimes D& 
4.\: =\:  \tilR_{(A| B\otimes C, D)}
\end{array}  \]
}

\vbox{
$((\bullet\otimes \bullet) \otimes (\bullet\otimes \bullet))$ \medskip

\[
\begin{tikzcd}[column sep = tiny]
    ABCD 
    \arrow[rrrr]
    \arrow[ddr, end anchor = {[xshift=1.7em]north west}]
    \arrow[dddrr, start anchor = {[xshift=-0.5em]south east}, end anchor = {[xshift=-0.2em]north}]
    \arrow[drr, dashed]
    &&&&
    CDAB
    \\
    &&
    ACDB
    \arrow[dr, dashed, end anchor = {[xshift=1.2em]north}]
    \arrow[urr, dashed]
    \arrow[ull, phantom, "3.", pos = 0.13, yshift=-1.7ex]
    \arrow[urr, phantom, "4.", pos = 0.13, yshift=-1.7ex]
    \arrow[dd, phantom, "5.", pos = -0.4]
    \arrow[dd, phantom, "6.", pos = 0.15]
    \\
    &
    ACBD\phantom{MMMM}
    \arrow[dr, start anchor = {[xshift=-1.2em]south}, end anchor = {[yshift=0.5ex]west}]
    \arrow[ur, dashed, start anchor = {[xshift=-1.2em]north}]
    &&
    \phantom{MMMM}CADB
    \arrow[uur, start anchor = {[xshift=-1.7em]north east}]
    \\
    &&
    CABD
    \arrow[ur, start anchor = {[yshift=0.5ex]east}, end anchor = {[xshift=1.2em]south}]
    \arrow[uuurr, start anchor = {[xshift=0.2em]north}, end anchor = {[xshift=0.5em]south west}]
    \arrow[uuull, phantom, "1.", pos = 0.1]
    \arrow[uuurr, phantom, "2.", pos = 0.1]
    \arrow[uu, phantom, "7.", pos = 0.2]
\end{tikzcd}
\]

\[ \begin{array}{lll}
1.\: =\:  \tilR_{(A,B|C)} \otimes D & 
2.\: =\:  C \otimes \tilR_{(A,B|D)}  &
3.\: =\:  A \otimes \tilR_{(B|C,D)} \\ 
4.\: =\:  \tilR_{(A| C, D)} \otimes B &
5.\: =\:  \tilR_{(A,B|C\otimes D)} &
6.\: =\: \otimes_{(R_{A,C},R_{B,D})}  \\
7.\: =\:  \tilR_{(A \otimes B|C,D)} &
\end{array} \]
}

Second, for any objects $A,B,C \in \C$, we define two 2-isomorphisms 
corresponding
to two proofs of the Yang--Baxter hexagon in a braided monoidal
category: \medskip

\[
\begin{tikzcd}[row sep = large, column sep = scriptsize]
    &
    BAC
    \arrow[r]
    \arrow[dd, phantom, pos=0.1, "\scriptstyle{\zwei \tilR_{(A|B,C)}^{-1}}"]
    \arrow[ddr, phantom, "\scriptstyle{\zwei R_{(A,R_{B,C})}^{-1}}"]
    &
    BCA
    \arrow[dr, start anchor = east]
    \arrow[dd, phantom, pos=0.9, "\scriptstyle{\zwei \tilR_{(A|B,C)}}"]
    \\
    ABC
    \arrow[ur, end anchor = west]
    \arrow[dr, end anchor = west]
    \arrow[urr, start anchor = east, end anchor = south]
    &&&
    CBA
    \\
    &
    ACB
    \arrow[r]
    \arrow[urr, start anchor = north, end anchor = west]
    &
    CAB
    \arrow[ur, start anchor = east]
\end{tikzcd}
\qquad
\begin{tikzcd}[row sep = large, column sep = scriptsize]
    &
    BAC
    \arrow[r]
    \arrow[drr, start anchor = south, end anchor = west]
    &
    BCA
    \arrow[dr, start anchor = east]
    \\
    ABC
    \arrow[ur, end anchor = west]
    \arrow[dr, end anchor = west]
    \arrow[drr, start anchor = east, end anchor = north]
    &&&
    CBA
    \\
    &
    ACB
    \arrow[uu, phantom, pos=0.1, "\scriptstyle{\zwei \tilR_{(A,B|C)}^{-1}}"]
    \arrow[uur, phantom, "\scriptstyle{\zwei R_{(R_{A,B},C)}}"]
    \arrow[r]
     &
    CAB
    \arrow[ur, start anchor = east]
    \arrow[uu, phantom, pos=0.9, "\scriptstyle{\zwei \tilR_{(A,B|C)}}"]
\end{tikzcd}
\]

\smallskip

We refer to these 2-morphisms as $S^+_{A,B,C}$ and $S^-_{A,B,C}$, respectively.
We require them to be equal:
\medskip
$(S^+ = S^-)$: \medskip

\[
\begin{tikzcd}[row sep = large, column sep = scriptsize]
    &
    BAC
    \arrow[r]
    \arrow[drr, start anchor = south]
    &
    BCA
    \arrow[dr, start anchor = east]
    \\
    ABC
    \arrow[ur, end anchor = west]
    \arrow[dr, end anchor = west]
    \arrow[drr, start anchor = east, end anchor = north]
    \arrow[urr, dashed, start anchor = east, end anchor = south]
    &&&
    CBA
    \\
    &
    ACB
    \arrow[r]
    \arrow[urr, dashed, start anchor = north, end anchor = west]
    &
    CAB
    \arrow[ur, start anchor = east]
\end{tikzcd}
\]

\smallskip

\end{defn}    

We can unpack this definition to obtain an explicit list of operations
and laws. In this form the definition is essentially due to Kapranov and 
Voevodsky, though with certain differences, which we list at the end of 
this section.

\begin{lma}
A braided monoidal 2-category $(\C,\otimes,1,R,\tilR_{(\stri|\stri,\stri)},
\tilR_{(\stri,\stri|\stri)})$ consists of the following data:

\begin{enumerate}

\item A semistrict monoidal 2-category $(\C,\otimes,1)$

\item $(\bullet \otimes \bullet)$ For any two objects
$A,B \in \C$ an equivalence $R_{A,B}\maps A \otimes B \to B \otimes A$
\item $({\to} \otimes \bullet)$ For any 1-morphism
$f: A \to A'$ and any object $B\in \C$ a 2-isomorphism 
\[
\begin{tikzcd}[row sep = huge, column sep = large]
    A \otimes B
    \arrow[r, "f \otimes B"]
    \arrow[d, swap, "R_{A,B}"]
    \arrow[dr, phantom, "\zwei R_{f,B}"]
    &
    A' \otimes B
    \arrow[d, "R_{A',B}"]
    \\
    B \otimes A
    \arrow[r, swap, "B \otimes f"]
    &
    B \otimes A'
\end{tikzcd}
\]

\item $(\bullet \otimes {\to})$ For any object $A\in \C$ and any 1-morphism
$ g\maps B \to B'$ a 2-isomorphism
\[
\begin{tikzcd}[row sep = huge, column sep = large]
    A \otimes B
    \arrow[r, "A \otimes g"]
    \arrow[d, swap, "R_{A,B}"]
    \arrow[dr, phantom, "\zwei R_{A,g}"]
    &
    A \otimes B'
    \arrow[d, "R_{A,B'}"]
    \\
    B \otimes A
    \arrow[r, swap, "g \otimes A"]
    &
    B' \otimes A
\end{tikzcd}
\]

\item $((\bullet \otimes \bullet) \otimes \bullet)$
For any objects $A,B,C \in \C$ a 2-iso
\[
\begin{tikzcd}[row sep = huge]
    A \otimes B \otimes C
    \arrow[rr, "R_{A, B \otimes C}"]
    \arrow[dr, swap, "R_{A,B} \otimes C"]
    & \phantom{2-iso} &
    B \otimes C \otimes A
    \\
    &
    B \otimes A \otimes C
    \arrow[u, phantom, yshift=1.4ex, "\Uparrow \tilR_{(A|B,C)}"]
    \arrow[ur, swap, "B \otimes R_{A,C}"]
\end{tikzcd}
\]

\item
$(\bullet \otimes (\bullet \otimes \bullet))$
For any objects $A,B,C \in \C$ a 2-isomorphism
\[
\begin{tikzcd}[row sep = huge]
    A \otimes B \otimes C
    \arrow[rr, "R_{A \otimes B, C}"]
    \arrow[dr, swap, "A \otimes R_{B,C}"]
    & \phantom{2-iso} &
    C \otimes A \otimes B
    \\
    &
    A \otimes C \otimes B
    \arrow[u, phantom, yshift=1.4ex, "\Uparrow \tilR_{(A,B|C)}"]
    \arrow[ur, swap, "R_{A,C} \otimes B"]
\end{tikzcd}
\]
\end{enumerate}

Moreover, these data must satisfy the following conditions:

$({\to} \otimes {\to})$ 
For any 1-morphisms $f \maps A \to A'$ and $g\maps B \to B'$ the following 
cube commutes:
\smallskip

\vbox{
\[
\begin{tikzcd}[row sep = large]
    &
    AB
    \arrow[dl]
    \arrow[rrr]
    \arrow[dd, dashed]
    & \phantom{1} & \phantom{6} &
    A'B
    \arrow[dl]
    \arrow[dd]
    \\
    AB'
    \arrow[rrr]
    \arrow[dd]
    &&&
    A'B'
    \arrow[dd]
    \arrow[ul, phantom, "1.", pos = 0.1]
    \arrow[dl, phantom, "5.", pos = 0.25]
    \\
    &
    BA
    \arrow[dl, dashed]
    \arrow[rrr, dashed]
    \arrow[dr, phantom, "2.", pos = 0.1]
    \arrow[uur, phantom, "6.", pos = 0.12]
    & \phantom{5} &&
    BA'
    \arrow[dl]
    \\
    B'A
    \arrow[rrr]
    \arrow[uuur, phantom, "3.", pos = 0.12]
    && \phantom{2} &
    B'A'
    \arrow[uuur, phantom, "4.", pos = 0.1]
\end{tikzcd}
\]
\bigskip

\[ 1.\: =\:  \otimes_{f,g}   \qquad
   2.\: =\:  \otimes_{g,f}  \qquad
   3.\: =\:    R_{A,g}          \qquad
   4.\: =\:   R_{A',g}             \qquad
   5.\: =\:  R_{f,B'}     \qquad
   6.\: =\:  R_{f,B}
\]
}
%

$(\bullet \otimes {\zwei})$ For any object $A \in \C$, any 1-morphisms
$f,f' \maps B \to B'$, and any 2-morphism
$\beta\maps f \Rightarrow f'$, the following prism commutes:
\smallskip

\[
\begin{tikzcd}[row sep = large, column sep = small]
    AB
    \arrow[ddr, phantom, yshift=1.5ex, "R_{A,f'}", pos = 0.3]
    \arrow[rr, bend left]
    \arrow[rr, bend right]
    \arrow[dd]
    &
    \zwei A \otimes \beta
    &
    AB'
    \arrow[dd]
    \\\\
    BA
    \arrow[rr, dashed, bend left]
    \arrow[rr, bend right]
    &
    \zwei \beta \otimes A
    &
    B'A
    \arrow[uul, phantom, yshift=-1.5ex, "R_{A,f}", pos = 0.3]
\end{tikzcd}\]
\smallskip

$({\zwei} \otimes \bullet)$
A similar prism, left to the reader. \smallskip

$(\to {\to} \otimes \bullet)$  For any pair of 1-morphisms
$ A \stackrel{f}{\to} A' \stackrel{f'}{\to} A''$ and any object $B \in \C$, 
the 2-isomorphism $R_{ff',B}$ coincides with the pasting \medskip

\[
\begin{tikzcd}[column sep = small]
    A \otimes B
    \arrow[rr]
    \arrow[dd]
    &&
    A' \otimes B
    \arrow[rr]
    \arrow[dd]
    \arrow[ddll, phantom, "\zwei R_{f,B}"]
    \arrow[ddrr, phantom, "\zwei R_{f',B}"]
    &&
    A'' \otimes B
    \arrow[dd]
    \\\\
    B \otimes A
    \arrow[rr]
    &&
    B \otimes A'
    \arrow[rr]
    &&
    B \otimes A''
\end{tikzcd}
\]
\medskip

$(\bullet \otimes {\to} \to)$ A similar pasting law, left to the reader.
\medskip 

$((\bullet \otimes \bullet) \otimes {\to})$
For any objects $A,B,C \in \C$ and any 1-morphism $f \maps C \to C'$, the 
following triangular prism commutes:

\smallskip
\vbox{
\[
\begin{tikzcd}[row sep = large]
    ABC
    \arrow [rrrr, "R_{AB,C}"]
    \arrow [dd, "AB \otimes f", swap]
    \arrow [drr]
    \arrow [dddrr, phantom, "\Uparrow 1.", pos = 0.14]
    \arrow [drrrr, phantom, yshift=0.5ex, "\Uparrow 3.", pos = 0.15]
    &&&&
    CAB
    \arrow [dd, "f \otimes AB"]
    \\
    &&
    ACB
    \arrow[urr]
    \arrow[dd]
    \arrow[drr, phantom, "\Uparrow 2.", pos = 0.1]
    && \phantom{3}
    \\
    ABC'
    \arrow[rrrr, dashed]
    \arrow[drr]
    \arrow[drrrr, phantom, yshift=0.5ex, "\Uparrow 4.", pos = 0.15]
    &&&&
    C'AB
    \arrow[ull, phantom, xshift=0.75em, "\Uparrow 5.", pos = 0.2]
    \\
    &&
    AC'B
    \arrow[urr]
    && \phantom{4}
\end{tikzcd}
\]
\bigskip

\[ 1.\: =\:  A \otimes R_{B,f} \qquad 
   2.\: =\:  R_{A,f} \otimes B  \qquad
   3.\: =\: \tilR_{(A,B|C)} \qquad
   4.\: =\: \tilR_{(A,B|C')} \qquad
   5.\: =\: R_{AB,f} \]
}

$({\to} \otimes (\bullet \otimes \bullet))$
A similar prism, left to the reader.

$(({\to} \otimes \bullet) \otimes \bullet)$
For any objects $A,B,C \in \C$ and any 1-morphism $f \maps A \to A'$, the 
following triangular prism commutes:

\medskip
\vbox{
\[
\begin{tikzcd}[row sep = large]
    ABC
    \arrow [rrrr, "R_{AB,C}"]
    \arrow [dd, "f \otimes BC", swap]
    \arrow [drr]
    \arrow [dddrr, phantom, "\Uparrow 1.", pos = 0.14]
    \arrow [drrrr, phantom, yshift=0.5ex, "\Uparrow 3.", pos = 0.15]
    &&&&
    CAB
    \arrow [dd, "C \otimes f \otimes B"]
    \\
    &&
    ACB
    \arrow[urr]
    \arrow[dd]
    \arrow[drr, phantom, "\Uparrow 2.", pos = 0.1]
    && \phantom{3}
    \\
    A'BC
    \arrow[rrrr, dashed]
    \arrow[drr]
    \arrow[drrrr, phantom, yshift=0.5ex, "\Uparrow 4.", pos = 0.15]
    &&&&
    CA'B
    \arrow[ull, phantom, xshift=0.75em, "\Uparrow 5.", pos = 0.2]
    \\
    &&
    A'CB
    \arrow[urr]
    && \phantom{4}
\end{tikzcd}
\]
\bigskip

\[ 1.\: =\:  \otimes_{(f,R_{B,C})}   \qquad
   2.\: =\:   R_{f,C} \otimes B  \qquad
   3.\: =\:  \tilR_{(A,B|C)}            \qquad
   4.\: =\: \tilR_{(A',B|C)}             \qquad
   5.\: =\:  R_{f\otimes B,C}
   \]
}

$((\bullet \otimes {\to}) \otimes \bullet)$,
$(\bullet \otimes ({\to} \otimes \bullet))$ and 
$(\bullet \otimes (\bullet \otimes {\to}))$
Similar prisms, left to the reader. \bigskip

$((\bullet \otimes \bullet \otimes \bullet) \otimes \bullet)$, 
$(\bullet \otimes (\bullet \otimes \bullet \otimes \bullet))$, 
$((\bullet \otimes \bullet) \otimes (\bullet \otimes \bullet))$ 
As in Definition \ref{bm2}.

\bigskip
$S^+ = S^-$ As in Definition \ref{bm2}.  

\end{lma}

\begin{proof} The 1-equivalences $R_{A,B}$ and 2-isomorphisms $R_{f,B}$ and
$R_{A,g}$ comprise the pseudonatural equivalence $R \maps \otimes \to
\otimes^\op$, and conditions $({\to} \otimes {\to})$, $({\bullet} \otimes
{\zwei})$, $({\zwei} \otimes {\bullet})$, $(\to {\to} \otimes {\bullet})$ and
$({\bullet} \otimes {\to} \to)$ state that it is indeed a pseudonatural
transformation.  The 2-morphisms $\tilR_{(A|B,C)}$ and $\tilR_{(A,B|C)}$
comprise the invertible modifications $\tilR_{(\stri|\stri ,\stri)}$ and
$\tilR_{(\stri,\stri|\stri)}$, and the commuting triangular prisms state
that these are indeed modifications, expressing naturality in each
argument.  The remaining 4 conditions come from Definition \ref{bm2}.
\end{proof}

Note that by $(\to {\to} \otimes \bullet)$ resp. $(\bullet \otimes {\to} \to)$ 
and by the invertibility of the respective 2-morphisms,
for any 
objects $A,B \in \C$ we have $R_{A,1_B} = 1_{R_{A,B}}$ and  $R_{1_A,B} =
1_{R_{A,B}}$.    

The above lemma makes it clear that our definition of braided monoidal
2-category differs from that of Kapranov and Voevodsky in precisely the
following points:
\begin{enumerate}
\item   Invertibility of the braiding.  Our definition
implies that the 1-morphisms $R_{A,B}$ are 
equivalences.   Kapranov and Voevodsky make no invertibility assumptions
on these 1-morphisms.  Our definition would agree with theirs on this
point, and otherwise stay the same, if we required $R \maps \otimes
\to \otimes^\op$ to be merely a pseudonatural transformation, rather than
a pseudonatural equivalence.  
\item  $S^+ = S^-$.  As already noted, Kapranov and Voevodsky omit
this condition. 
\item  Naturality of $\tilR_{(\stri|\stri ,\stri)}$ and $\tilR_{(\stri,\stri | \stri)}$.  
Our definition implies the commutativity of 6 triangular prisms
expressing the naturality in each argument of these modifications.
Kapranov and Voevodsky substitute cubes for 4 of these prisms, namely
$(\bullet \otimes ({\to} \otimes \bullet))$, $(\bullet \otimes (\bullet \otimes {\to}))$, 
$((\bullet \otimes {\to}) \otimes \bullet)$ and $(({\to} \otimes \bullet) \otimes \bullet)$.
By the following lemma one can deduce these cubes from 
the remaining data
--- but not, it appears, vice versa.  In personal communication,
Kapranov agreed that all these prisms should hold.  
\end{enumerate}

\begin{lma}      \label{Lcub}
For any three objects $A,B,C \in \C$ and any morphism $f\maps B \to B'$,
the following cube commutes. \bigskip

\vbox{
\[
\begin{tikzcd}[row sep = large]
    &
    ABC
    \arrow[dl]
    \arrow[rrr]
    \arrow[dd, dashed]
    & \phantom{1} & \phantom{6} &
    AB'C
    \arrow[dl]
    \arrow[dd]
    \\
    ACB
    \arrow[rrr]
    \arrow[dd]
    &&&
    ACB'
    \arrow[dd]
    \arrow[ul, phantom, "\scriptstyle{1.}", pos = 0.1]
    \arrow[dl, phantom, "\scriptstyle{5.}", pos = 0.25]
    \\
    &
    BCA
    \arrow[dl, dashed]
    \arrow[rrr, dashed]
    \arrow[dr, phantom, "\scriptstyle{2.}", pos = 0.1]
    \arrow[uur, phantom, "\scriptstyle{6.}", pos = 0.12]
    & \phantom{5} &&
    B'CA
    \arrow[dl]
    \\
    CBA
    \arrow[rrr]
    \arrow[uuur, phantom, "\scriptstyle{4.}", pos = 0.12]
    && \phantom{2} &
    CB'A
    \arrow[uuur, phantom, "\scriptstyle{3.}", pos = 0.1]
\end{tikzcd}
\]
\smallskip

\[ \begin{array}{lll}
   1.\: =\:  A \otimes R_{f,C}   &
   2.\: =\:  R_{f,C} \otimes A   &
   3.\: =\:    R_{A,R_{B',C}}  \\
   4.\: =\:   R_{A,R_{B,C}}   &          
   5.\: =\:  R_{A,C \otimes f}    & 
   6.\: =\:  R_{A,f \otimes C}
\end{array}\]
}
\end{lma}

\begin{proof} 
This is an special case of the axiom $(\bullet \otimes {\zwei})$ together with $(\bullet \otimes {\to} \to )$.   
\end{proof}

We refer to this cube with the hieroglyph 
$(\bullet \otimes ({\to} \otimes \bullet))'$. One can similarly
prove the analogous cube corresponding to the hieroglyph 
$(\bullet \otimes (\bullet \otimes {\to}))'$ commutes.  
Moreover, we can prove the commutativity of cubes
corresponding to the hieroglyphs 
$((\bullet \otimes {\to}) \otimes \bullet)'$ and 
$(({\to} \otimes \bullet) \otimes \bullet)'$ using 
$({\zwei} \otimes \bullet)$  and  $(\to {\to} \otimes \bullet)$.

%% file: 3centerconstruction.tex
\section{The Center Construction}

Let $(\C,\otimes,1)$ be a semistrict monoidal 2-category.  The center $\Z(\C)$ would be easy to construct if we had a properly functioning theory of semistrict weak 4-categories. As it stands, all we can do is use our limited insight into 4-categories to guess the right answer, and then try to justify it by proving that we obtain a braided monoidal 2-category with good properties. We proceed in several stages.  First we describe  $\Z(\C)$ as a 2-category. Then we describe the monoidal structure, and then  the braiding.   

\subsection{$\Z(\C)$ as a 2-Category}

As noted in Section~\ref{center}, the center construction 
applied to a monoid yields its usual center, because a certain 
square must commute.    However, as one would expect from the weakening
principle, when 
$\C$ is a monoidal category the corresponding square need only commute \emph{up
to a specified natural isomorphism}.  An object of $\Z(\C)$ thus turns out to be
an object $A \in \C$ equipped with a natural isomorphism $R_{A, \stri} 
\maps A \otimes \stri \To \stri \otimes A$ satisfying various coherence laws,
such as the commutativity of following diagram:
\[
\begin{tikzcd}[row sep = large]
    A \otimes X \otimes Y
      \arrow[rr, "R_{A,X \otimes Y}"]
      \arrow[ddr, swap, "R_{A,X} \otimes Y"]
    &&
    X \otimes Y \otimes A
    \\\\ &
    X \otimes A \otimes Y
      \arrow[uur, swap, "X \otimes R_{A,Y}"]
\end{tikzcd}
\]
for any objects $X,Y \in \C$.  Of course, this diagram is part of
the definition of a braided monoidal category.  Similarly, the morphisms 
in $\Z(\C)$ also work out to have properties that form part of the definition
of a braided monoidal category.  

Heuristic 4-categorical computations suggest how these
patterns should continue when $\C$ is a monoidal 2-category.  We thus define
$\Z(\C)$ as follows.  

\medskip
\define{Objects in $\Z(\C)$:}
\medskip

An object of $\Z(\C)$ is a triple $(A,R_{A,\stri},
\tilR_{(A|\stri,\stri)})$ consisting of:
\begin{enumerate}
\item
an object $A \in \C$
\item
a pseudonatural equivalence 
$ R_{A,\stri} \maps A \otimes \stri \To \stri\otimes A$
\item
an invertible modification $\tilR_{(A|\stri,\stri)}$, giving for any objects
$X,Y \in \C$ a 2-isomorphism

\[
\begin{tikzcd}[row sep = large]
    A \otimes X \otimes Y
      \arrow[rr, "R_{A,X \otimes Y}"]
      \arrow[ddr, swap, "R_{A,X} \otimes Y"]
    & \phantom{R} &
    X \otimes Y \otimes A
    \\\\ &
    X \otimes A \otimes Y
      \arrow[uur, swap, "X \otimes R_{A,Y}"]
      \arrow[uu, phantom, "\Uparrow \tilR_{(A|X,Y)}", pos = 0.65]
\end{tikzcd}
\]
\end{enumerate}
such that for any objects $X,Y,Z \in \C$, the tetrahedron
$(\bullet \otimes(\bullet \otimes \bullet \otimes \bullet))$ commutes.
\medskip

Here we mean that the diagram 
$(\bullet \otimes(\bullet \otimes \bullet \otimes \bullet))$ 
commutes with objects $A,X,Y,Z$,
and with the modification $\tilR_{(\stri|\stri,\stri)}$ in the 
definition of a braided monoidal 2-category replaced by the above
$ \tilR_{(A|\stri,\stri)} $.
Throughout the following we use the hieroglyphical notation in this way.  
Also, we use letters near the beginning of the alphabet to denote objects of $\C$ 
underlying objects in $\Z(\C)$, and letters near the end to denote objects
of $\C$ being used as such.  

\begin{bem}  \label{objects} 
The fact that $R_{A,\stri}$ is a pseudonatural equivalence can be 
expressed equivalently as follows: 
for any object $X \in \C$, there exists an equivalence 
$R_{A,X} \maps A \otimes X \to  X \otimes A$,
and for any morphism $f \maps X \to Y$ in $\C$, there exists a 2-isomorphism
$R_{A,f} \maps (A \otimes f) \circ R_{A,Y} \To R_{A,X}\circ (f \otimes A) $:

\[
\begin{tikzcd}[row sep = 4.5 em, column sep = 3.2 em]
    A \otimes X
      \arrow[r, "R_{A,X}"]
      \arrow[d, swap, "A \otimes f"]
      \arrow[dr, phantom, "\Uparrow R_{A,f}"]
    &
    X \otimes A
      \arrow[d, "f \otimes A"]
    \\
    A \otimes Y
      \arrow[r, swap, "R_{A,Y}"]
    &
    Y \otimes A
\end{tikzcd}
\]
such that 
$(\bullet \otimes {\to} \to)$ and $(\bullet \otimes {\zwei})$ commute.
\smallskip

Similarly, the fact that  
$\: \tilR_{(A|\stri , \stri)}$ is a modification means that the diagrams
$(\bullet \otimes ({\to} \otimes \bullet))$ and
$(\bullet \otimes (\bullet \otimes {\to}))$ commute.
\end{bem}

\define{Morphisms in $\Z(\C)$:}
\medskip

A morphism in $\Z(\C)$ from
$(A,R_{A,\stri},\tilR_{(A|\stri,\stri)})$ to $(B,R_{B,\stri},\tilR_{(B|\stri,\stri)})$ 
is a pair $(f, R_{f, \stri})$ consisting of:
\begin{enumerate}
\item
a morphism $f\maps A \to B$
\item
an invertible modification $R_{f, \stri}$, giving for any object $X \in \C$ a
2-isomorphism

\[
\begin{tikzcd}[row sep = 4.5 em, column sep = 3.2 em]
    A \otimes X
      \arrow[r, "f \otimes X"]
      \arrow[d, swap, "R_{A,X}"]
      \arrow[dr, phantom, "\zwei R_{f,X}"]
    &
    B \otimes X
      \arrow[d, "R_{B,X}"]
    \\
    X \otimes A
      \arrow[r, swap, "X \otimes f"]
    &
    X \otimes B
\end{tikzcd}
\]

\end{enumerate}
such that the prism $({\to} \otimes (\bullet \otimes \bullet))$ commutes.

\begin{bem}
The fact that $R_{f,\stri}$ is a modification can be expressed equivalently 
by saying that $({\to} \otimes {\to})$ commutes. (Note that 
$(f \otimes \stri)R_{B,\stri}$ and $R_{A,\stri}(\stri \otimes f)$ are 
pseudonatural transformations in an obvious way.) 
\end{bem}

\medskip
\define{2-Morphisms in $\Z(\C)$:}
\medskip

A 2-morphism $\alpha$ in $\Z(\C)$ from 
$(f,R_{f,\stri})$ to $(g,R_{g,\stri})$ is 
\begin{enumerate}
\item a 2-morphism $\alpha \maps f \To g$ in $\C$ 
\end{enumerate}
such that $({\zwei} \otimes \bullet)$ commutes.
\bigskip

We define the composition operations in $\Z(\C)$ as follows.  
Composition of morphisms is defined by:
\[ (f,R_{f,\stri}) \circ (g,R_{g,\stri}) := (f \circ g, ((f\otimes \stri) \circ R_{g,\stri})
\cdot (R_{f,\stri} \circ (\stri \otimes g)) ) \]
where $f\maps A \to B$ and $g \maps B \to C$ are the underlying 1-morphisms 
in $\C$.   Note that for any object $X \in \C$, the 2-morphism 
$((f\otimes X) \circ R_{g,X}) \cdot (R_{f,X} \circ (X \otimes g))$ 
equals the back of the following diagram:

\medskip
\[
\begin{tikzcd}[row sep = 2.3 em, column sep = 2.7 em]
    AX
    \arrow [rrrr, "(f \circ g) \otimes X"]
    \arrow [dd, swap, "R_{A,X}"]
    \arrow [drr, swap, "f \otimes X"]
    \arrow [dddrr, phantom, "\zwei R_{f,X}"]
    \arrow [drrrr, phantom, yshift=0.5ex, "\Uparrow \textrm{id}", pos = 0.15]
    &&&&
    CX
    \arrow [dd, "R_{C,X}"]
    \\ &&
    BX
    \arrow[urr, swap, "g \otimes X"]
    \arrow[dd]
    \arrow[drr, phantom, "\zwei R_{g,X}"]
    && \phantom{3}
    \\
    XA
    \arrow[rrrr, dashed]
    \arrow[drr, swap, "X \otimes f"]
    \arrow[drrrr, phantom, yshift=0.5ex, "\Uparrow \textrm{id}", pos = 0.15]
    &&&&
    XC
    \\ &&
    XB
    \arrow[urr, swap, "X \otimes g"]
    && \phantom{4}
\end{tikzcd}
\]
\bigskip

\begin{bem}
Eventually 
this will imply that the braiding in $\Z(\C)$ satisfies $(\to {\to} \otimes \bullet)$.
\end{bem}

To show that the composite of morphisms in $\Z(\C)$ is again a morphism,
we have to check that $({\to} \otimes {\to})$ and 
$({\to} \otimes (\bullet \otimes \bullet))$ commute.
These can be seen by pasting together two diagrams of the 
form $ ({\to} \otimes {\to})$ and $({\to} \otimes (\bullet \otimes \bullet))$, 
respectively.

Vertical and horizontal composition of 2-morphisms is defined the same 
as in $\C$; one can check that these composites again satisfy 
$({\zwei} \otimes \bullet)$ by pasting together two diagrams of this form.  

\subsection{The Monoidal Structure}
\label{monoidal}

We have to show that $\Z(\C)$ bears a monoidal structure $(\Z(\C),
\otimes_{\Z(\C)}, I) $, such that all the requirements for a monoidal category 
given in Definition \ref{mon} are satisfied.
\bigskip

(Ad 4.1): The object $I \in \Z(\C)$ is 
$(I,1_{\stri},1_{1_{(\stri \otimes \stri)}})$.

\bigskip
\define{The tensor product of objects:}
(Ad 4.2): The tensor product of two objects
$(A,R_{A,\stri},\tilR_{(A|\stri ,\stri)}) \otimes_{\Z(\C)} 
(B,R_{B,\stri},\tilR_{(B|\stri, \stri)}) $ 
is defined to be the triple 
$(A \otimes B,(R_A \otimes R_B)_{\stri}, 
(\tilR_A \otimes \tilR_B)_{(\stri, \stri)})$, where:
\begin{enumerate}
\item The underlying $\C$-object is the tensor product $A \otimes B$ in $\C$.
\item By Remark~(\ref{objects}), the underlying pseudonatural equivalence
$(R_A \otimes R_B)_{\stri}: (A \otimes B) \otimes \stri 
\To \stri \otimes (A \otimes B) $ assigns a 1-morphism $(R_A \otimes R_B)_X$
to any object $X \in \C$ and a 2-morphism $(R_A \otimes R_B)_f$ to any 
1-morphism $f \maps X \to Y$.  These are given as follows:
\[    (R_A \otimes R_B)_X  = (A \otimes R_{B,X})(R_{A,X} \otimes B), \]
\[(R_A \otimes R_B)_f = ((A \otimes R_{B,f}) \circ (R_{A,Y}\otimes B))
\cdot (A \otimes R_{B,X})\circ (R_{A,f} \otimes B)) ,\]
or in terms of a diagram:
\[
\begin{tikzcd}[row sep = huge, column sep = large]
    ABX
      \arrow[d, swap, "AB \otimes f"]
      \arrow[r, "A \otimes R_{B,X}"]
      \arrow[dr, phantom, "\Uparrow A \otimes R_{B,f}"]
    &
    AXB
      \arrow[d]
      \arrow[r, "R_{A,X} \otimes B"]
      \arrow[dr, phantom, "\Uparrow R_{A,f} \otimes B"]
    &
    XAB
      \arrow[d, "XA \otimes f"]
    \\
    ABY
      \arrow[r, swap, "A \otimes R_{B,Y}"]
    &
    AYB
      \arrow[r, swap, "R_{A,Y} \otimes B"]
    &
    YAB
\end{tikzcd}
\]
\smallskip

\begin{bem} This will imply that the braiding in $\Z(\C)$ satisfies
 $((\bullet \otimes \bullet) \otimes {\to})$.   \end{bem}

To show that these data constitute a pseudonatural equivalence, 
we have to show that 
$(\bullet \otimes {\to} \to)$ and $(\bullet \otimes {\zwei})$ hold.
This can be done easily by pasting together the corresponding diagrams for
$R_{A,\stri}$ and $R_{B,\stri}$.
 
\item The underlying modification 
$(\tilR_A \otimes \tilR_B)_{(X,Y)}: (A \otimes R_{B,X} \otimes Y) 
(R_{A,X} \otimes B \otimes Y)
(X \otimes A \otimes R_{B,Y})(X \otimes R_{A,Y} \otimes B)  \To
(A \otimes R_{B,X \otimes Y}) (R_{A,X \otimes Y} \otimes B)$ 
is defined to be the pasting:

\[
\begin{tikzcd}[column sep = 0.9 em]
    ABXY
      \arrow[drrr, "A \otimes R_{B,XY}"]
      \arrow[ddr]
      \arrow[dddrrrr, phantom, "\Uparrow A \otimes \tilR_{(B|X,Y)}", pos = 0.31]
    \\
    &&&
    AXYB
      \arrow[drrr, "R_{A,XY} \otimes B"]
      \arrow[ddr]
    \\
    &
    AXBY
      \arrow[urr]
      \arrow[ddr]
      \arrow[drrr, phantom, "\Uparrow \otimes_{R_{A,X},R_{B,Y}}^{-1}", pos = 0.45]
    &&&& \phantom{ABXY} &
    XYAB
      \arrow[lllll, phantom, "\Uparrow \tilR_{(A|X,Y)} \otimes B", pos = 0.31]
    \\
    &&&&
    XAYB
      \arrow[urr]
    \\
    &&
    XABY
      \arrow[urr]
\end{tikzcd}
\]
\bigskip

Again it is easy to verify that this satisfies 
$(\bullet \otimes ({\to} \otimes \bullet))$ and 
$(\bullet \otimes (\bullet \otimes {\to}))$ and hence is a modification.
\end{enumerate}

To show that this definition gives in fact an object in $\Z(\C)$, we have to 
verify that 
$(\bullet \otimes (\bullet \otimes \bullet \otimes \bullet))$ is satisfied.
The following picture shows the tetrahedron. 
Those vertices in the picture that are vertices of the tetrahedron are 
written in big capitals. The remaining vertices occur since they are needed 
for the decomposition.
\medskip

\[
\begin{tikzcd}[row sep = small]
    &&
    XYZAB
    \\ \phantom{A} \\ &
    \scriptstyle{AXYZB}
      \arrow[uur]
    &&
    \scriptstyle{XYAZB}
      \arrow[uul]
    \\ \phantom{A} &&
    \scriptstyle{XAYZB}
      \arrow[uuu]
    \\
    ABXYZ
      \arrow[uur]
      \arrow[dr]
      \arrow[r, dashed]
    &
    \scriptstyle{AXYBZ}
      \arrow[rrr, dashed]
    &&&
    XYABZ
      \arrow[uul]
    \\ &
    \scriptstyle{AXBYZ}
      \arrow[dr]
    &&
    \scriptstyle{XAYBZ}
      \arrow[ur]
    \\ &&
    XABYZ
      \arrow[ur]
      \arrow[uuu]
\end{tikzcd}
\]

The following picture gives a decomposition of the
tetrahedron into four smaller commutative diagrams. \smallskip
\begin{center}
\begin{tikzpicture}
    \node (1a) at (6,9.5)    {$AXYZB$};
    \node (1b) at (9,8.5)    {$XAYZB$};
    \node (1c) at (11.5,9.5) {$XYAZB$};
    \node (1d) at (9.5,11)   {$XYZAB$};

    \node (2a) at (0,6)      {$ABXYZ$};
    \node (2b) at (2,4)      {$AXBYZ$};
    \node (2c) at (4,5)      {$AXYBZ$};
    \node (2d) at (3,7)      {$AXYZB$};

    \node (3a) at (4.5,3)    {$AXBYZ$};
    \node (3b) at (7.5,2)    {$XABYZ$};
    \node (3c) at (6.5,4)    {$AXYBZ$};
    \node (3bc) at (9.5,3)   {$XAYBZ$};
    \node (3d) at (5.5,6)    {$AXYZB$};
    \node (3e) at (8.5,5)    {$XAYZB$};

    \node (4a) at (7.8,6)    {$AXYBZ$};
    \node (4b) at (10.8,5)   {$XAYBZ$};
    \node (4c) at (6.8,8)    {$AXYZB$};
    \node (4bc) at (9.8,7)   {$XAYZB$};
    \node (4d) at (13.3,6)   {$XYABZ$};
    \node (4e) at (12.3,8)   {$XYAZB$};

    \draw[->]        (1a)  to (1b);
    \draw[dashed,->] (1a)  to (1c);
    \draw[->]        (1a)  to (1d);
    \draw[->]        (1b)  to (1c);
    \draw[->]        (1b)  to (1d);
    \draw[->]        (1c)  to (1d);

    \draw[->]        (2a)  to (2b);
    \draw[dashed,->] (2a)  to (2c);
    \draw[->]        (2a)  to (2d);
    \draw[->]        (2b)  to (2c);
    \draw[->]        (2b)  to (2d);
    \draw[->]        (2c)  to (2d);

    \draw[->]        (3a)  to (3b);
    \draw[dashed,->] (3a)  to (3c);
    \draw[->]        (3a)  to (3d);
    \draw[->]        (3b)  to (3bc);
    \draw[->]        (3b)  to (3e);
    \draw[dashed,->] (3c)  to (3bc);
    \draw[dashed,->] (3c)  to (3d);
    \draw[->]        (3bc) to (3e);
    \draw[->]        (3d)  to (3e);

    \draw[->]        (4a)  to (4b);
    \draw[->]        (4a)  to (4c);
    \draw[dashed,->] (4a)  to (4d);
    \draw[->]        (4b)  to (4bc);
    \draw[->]        (4b)  to (4d);
    \draw[->]        (4c)  to (4bc);
    \draw[->]        (4c)  to (4e);
    \draw[->]        (4bc) to (4e);
    \draw[->]        (4d)  to (4e);
\end{tikzpicture}
\end{center}

 \medskip

Two of them are tetrahedra of the form
$(\bullet \otimes (\bullet \otimes \bullet \otimes \bullet))$, 
tensored by an object from the left and the right, respectively.
The upper of the two triangular prisms commutes by the axioms $4.(vii)$ 
with $\alpha = \tilR_{(A|X,Y)}$ and $g = R_{Z,B}$, 
together with $4.(viii)$. 
The lower commutes by $4.(vi)$, applied to $\beta = \tilR_{(B|Y,Z)}$ and 
$f = R_{A,X}$ together with $4.(viii)$. 
One can verify that this tensor product is in fact associative.  

We shall often write $(A\otimes B, R_A \otimes R_B , \tilR_A \otimes \tilR_B)$
as a shorthand symbol for the tensor product of objects in $\Z(\C)$.
\bigskip

\define{The tensor product of an object and a morphism:}

(Ad 4.3):
Let $(f,R_{f,\stri}) \maps (A,R_{A,\stri},\tilR_{(A|\stri,\stri)}) 
\to (A',R_{A',\stri},\tilR_{(A'|\stri,\stri)})$ 
be a morphism in $\Z(\C)$ and let $(B,R_{B,\stri},\tilR_{(B|\stri,\stri)})$ be an object.
Their tensor product is the morphism given by the pair 
\[(f \otimes B,(\otimes _{(f,R_{B,\stri})}\circ (R_{A',\stri} \otimes B)) 
\cdot((A \otimes R_{B,\stri}) \circ (R_{f,\stri} \otimes B))),\]
or in terms of a diagram:

\[
\begin{tikzcd}[row sep = huge, column sep = large]
    ABX
      \arrow[d, swap, "f \otimes BX"]
      \arrow[r, "A \otimes R_{B,X}"]
      \arrow[dr, phantom, "\Uparrow \otimes_{f,R_{B,X}}"]
    &
    AXB
      \arrow[d]
      \arrow[r, "R_{A,X} \otimes B"]
      \arrow[dr, phantom, "\Uparrow R_{f,X} \otimes B"]
    &
    XAB
      \arrow[d, "X \otimes f \otimes B"]
    \\
    A'BX
      \arrow[r, swap, "A' \otimes R_{B,X}"]
    &
    A'XB
      \arrow[r, swap, "R_{A',X} \otimes B"]
    &
    XA'B
\end{tikzcd}
\]

\begin{bem}
This will imply that the braiding in 
$\Z(\C)$ satisfies $(({\to} \otimes \bullet) \otimes \bullet)$.  \end{bem}

\medskip

(Ad 4.4)
Let $(f,R_{f,\stri}) : (B,R_{B,\stri},\tilR_{(B|\stri,\stri)}) 
\to (B',R_{B',\stri},\tilR_{(B'|\stri,\stri)})$ 
be a morphism in $\Z(\C)$ and let $(A,R_{A,\stri},\tilR_{(A|\stri,\stri)})$ 
be an object.  Their tensor product is the pair
\[(A \otimes f,
((A \otimes R_{f,\stri}) \circ (R_{A,\stri} \otimes B'))  \cdot
((A \otimes R_{B,\stri}) \circ \otimes _{R_{A,\stri},f})),  \]
or in terms of a diagram:
\[
\begin{tikzcd}[row sep = huge, column sep = large]
    ABX
      \arrow[d, swap, "A \otimes f \otimes X"]
      \arrow[r, "A R_{B,X}"]
      \arrow[dr, phantom, "\Uparrow A \otimes R_{f,X}"]
    &
    AXB
      \arrow[d]
      \arrow[r, "R_{A,X}  B"]
      \arrow[dr, phantom, "\Uparrow \otimes_{R_{A,X},f}"]
    &
    XAB
      \arrow[d, "XA \otimes f"]
    \\
    AB'X
      \arrow[r, swap, "A R_{B',X}"]
    &
    AXB'
      \arrow[r, swap, "R_{A,X} B'"]
    &
    XAB'
\end{tikzcd}
\]

\begin{bem} This will imply that
the braiding in $\Z(\C)$ satisfies $((\bullet \otimes {\to}) \otimes \bullet)$.
\end{bem}

To verify that these formulas really define morphisms in $\Z(\C)$,
one must check that $({\to} \otimes {\to})$ and
$({\to} \otimes(\bullet \otimes \bullet))$ hold. 
We only do this for $4.4$; the other case being 
similar.   To show $({\to} \otimes {\to})$ one pastes together two cubes, 
one being the $({\to} \otimes {\to})$ cube for $f \maps B \to B'$ and
$g \maps X \to Y$
tensored on the left by $A$, the other being a special case of $5.(vii)$.
For $({\to} \otimes(\bullet \otimes \bullet))$
we must show the following diagram commutes:

\vbox{
\[
\begin{tikzcd}[column sep = 0.8 em, row sep = 1.0 em]
    ABXY
      \arrow[drrr, "A \otimes R_{B,XY}"]
      \arrow[ddr]
      \arrow[dddd, swap, "R_{A,B} \otimes XY"]
      \arrow[dddrrrr, phantom, "\Uparrow A \otimes \tilR_{(B|X,Y)}", pos = 0.31]%
      \arrow[ddddddddrr, phantom, "1." pos = 0.05]
    \\
    &&&
    AXYB
      \arrow[drrr, "R_{A,XY} \otimes B"]
      \arrow[ddr]
      \arrow[dddd, dashed]
      \arrow[dddddll, phantom, "7.", pos = 0.95]
      \arrow[ddddddr, phantom, "8.", pos = 0.92]
    \\
    &
    AXBY
      \arrow[urr]
      \arrow[ddr]
      \arrow[dddd]
      \arrow[drrr, phantom, "\Uparrow \otimes_{R_{A,X},R_{B,Y}}^{-1}", pos = 0.45]
      \arrow[ddddddr, phantom, "2.", pos = 0.06, xshift=0.3em]
    &&&& \phantom{ABXY} &
    XYAB
      \arrow[dddd, "XY \otimes R_{A,B}"]
      \arrow[lllll, phantom, "\Uparrow \tilR_{(A|X,Y)} \otimes B", pos = 0.31]
    \\
    &&&&
    XAYB
      \arrow[urr]
      \arrow[dddd]
      \arrow[dddrr, phantom, "4.", pos = 0.05]
    \\
    AB'XY
      \arrow[drrr, dashed]
      \arrow[ddr, swap, "A \otimes R_{B',X} \otimes Y"]
    &&
    XABY
      \arrow[urr]
      \arrow[dddd]
    \\
    &&&
    AXYB'
      \arrow[ddr, dashed]
      \arrow[drrr, dashed]
      \arrow[uuull, phantom, "5.", pos = 0.05]
    \\
    &
    AXB'Y
      \arrow[urr, dashed]
      \arrow[ddr, swap, "R_{A,X} \otimes B'Y"]
    &&&&&
    XYAB'
      \arrow[uuull, phantom, "6.", pos = 0.05]
    \\
    &&&&
    XAYB'
      \arrow[urr, swap, "X \otimes R_{A,Y} \otimes B'"]
    \\
    &&
    XAB'Y
      \arrow[urr, swap, "XA \otimes R_{B',Y}"]
      \arrow[uuur, phantom, "3.", pos = 0.1]
\end{tikzcd}
\]

\bigskip

\[ \begin{array}{llll}
1.\: =\: A \otimes R_{f,X} \otimes Y & 
2.\: =\:  \otimes_{R_{A,X},f\otimes Y} & 
3.\: =\: X \otimes A \otimes R_{f,Y}  &
4.\: =\: X \otimes \otimes_{R_{A,Y},f}  \\ 
5.\: =\: A \otimes R_{f,X \otimes Y} & 
6.\: =\: \otimes_{R_{A,X \otimes Y},f} &
7.\: =\: A \otimes X \otimes R_{f,Y}  &
8.\: =\:  \otimes_{R_{A,X},Y \otimes f}  
\end{array}    \]
}

We cut it into one rectangular and two triangular prisms.
To see that the left triangular prism commutes, we apply 
$({\to} \otimes (\bullet \otimes \bullet))$ to $(f,R_{f,\stri})$,
tensored on the left by $A$.
The rectangular prism commmutes by $5.(vi)$ and $5.(viii)$, applied to the 
2-morphism $R_{f,Y}$.
The right triangular prism commutes by $5.(iv)$ and $5.(vii)$, applied to the 
2-morphism $\tilR_{(A|X,Y)} $.

\bigskip

\define{The tensor product of an object and a 2-morphism:}

(Ad 4.5):
For any object $(A,R_{A,\stri},\tilR_{(A|\stri , \stri)})$ and any 
2-morphism 
$\alpha\maps (f,R_{f,\stri}) \To (f',R_{f',\stri})$ we have a 2-morphism
\[ A \otimes \alpha : (A \otimes f,\dots) \To (A \otimes f',\dots)\]

(Ad 4.6):
For any object $(B,R_{B,\stri},\tilR_{(B|\stri, \stri)})$ and any 2-morphism 
$\alpha: (g,R_{g,\stri}) \To (g',R_{g',\stri})$ we have a 2-morphism
\[ \alpha \otimes B : (g \otimes B,\dots) \To (g' \otimes B,\dots)\]
We must verify that these are 2-morphisms in $\Z(\C)$, so we must check
$({\zwei} \otimes \bullet)$. 
We do this only for $4.5$.

\bigskip

\vbox{
\[
\begin{tikzcd}[row sep = large, column sep = tiny]
    ABX
    \arrow[ddr, phantom, yshift=-1.5ex, "2.", pos = 0.1]
    \arrow[rr, bend left]
    \arrow[rr, bend right]
    \arrow[dd, swap, "A \otimes R_{B,X}"]
    &
    \zwei A \otimes \alpha \otimes X
    &
    AB'X
    \arrow[dd, "A \otimes R_{B',X}"]
    \\\\
    AXB
    \arrow[rr, dashed, bend left]
    \arrow[rr, bend right]
    \arrow[dd, swap, "R_{A,X} \otimes B"]
    \arrow[ddr, phantom, yshift=-1.5ex, "4.", pos = 0.1]
    &
    \zwei AX \otimes \alpha
    &
    AXB'
    \arrow[dd, "R_{A,X} \otimes B'"]
    \arrow[uul, phantom, yshift=1.5ex, "1.", pos = 0.1]
    \\\\
    XAB
    \arrow[rr, dashed, bend left]
    \arrow[rr, bend right]
    &
    \zwei XA \otimes \alpha
    &
    XAB'
    \arrow[uul, phantom, yshift=1.5ex, "3.", pos = 0.1]
\end{tikzcd}
\]

\[  \begin{array}{ll}
1.\: =\: A \otimes R_{f,X} & 2.\: =\: A \otimes R_{g,X} \\
3.\: =\: \otimes_{R_{A,X},f} & 4.\: =\:  \otimes_{R_{A,X},g}
\end{array}
\]
}

The upper prism commutes by $({\zwei} \otimes \bullet)$ tensored from the left 
by $A$. The lower prism commutes by an application of the axiom $4.(vi)$ for
monoidal $2$-categories to the 2-morphism $\alpha$.

\define{The tensor product of morphisms:}  

(Ad 4.7):
For any morphisms 
$(f,R_{f,\stri}) : (A,R_A,\tilR_A) \to (A',R_{A'},\tilR_{A'})$ and
$(g,R_{g,\stri}):(B,R_B,\tilR_B) \to (B',R_{B'},\tilR_{B'})$ 
we have a 2-isomorphism:
\[ \otimes_{(f,R_{f,\stri}),(g,R_{g,\stri})} := \otimes_{f,g} \] 
To verify that this is a 2-morphism in $\Z(\C)$, we have to check
$({\zwei} \otimes \bullet)$.
The following diagram gives the proof.

\vbox{
\[
\begin{tikzcd}[column sep = small, row sep = 1.6 em]
    ABX
      \arrow[ddd, swap, "A \otimes R_{B,X}"]
      \arrow[rrrr, "f \otimes BX"]
      \arrow[ddr]
      \arrow[ddrrrrr, phantom, "\Uparrow \otimes_{f,g} \otimes X", pos = 0.17]
      \arrow[dddddr, phantom, "1.", pos = 0.1]
    &&& \phantom{ABX} &
    A'BX
      \arrow[ddr, "A' \otimes g \times X"]
      \arrow[ddd, dashed]
      \arrow[dddddr, phantom, "3.", pos = 0.1]
    \\\\ &
    AB'X
      \arrow[ddd]
      \arrow[rrrr, "f \otimes B'X"]
      \arrow[dddrrrr, phantom, "2.", pos = 0.1, xshift=-1em]
    &&&&
    A'B'X
      \arrow[ddd, "A' \otimes R_{B',X}"]
    \\
    AXB
      \arrow[ddd, swap, "R_{A,X} \otimes B"]
      \arrow[ddr]
      \arrow[rrrr, dashed]
      \arrow[dddddr, phantom, "5.", pos = 0.1]
      \arrow[ddrrrrr, phantom, "\Uparrow \otimes_{f \otimes X, g}"]
    &&&&
    A'XB
      \arrow[ddd, dashed]
      \arrow[ddr, dashed]
      \arrow[uuullll, phantom, "4.", pos = 0.1, xshift=1.5em]
      \arrow[dddddr, phantom, "7.", pos = 0.1]
    \\\\ &
    AXB'
      \arrow[ddd]
      \arrow[rrrr, swap, "f \otimes XB'"]
      \arrow[dddrrrr, phantom, "6.", pos = 0.1, xshift=-1em]
    &&&&
    A'XB'
      \arrow[ddd, "R_{A',X} \otimes B'"]
    \\
    XAB
      \arrow[ddr, swap, "XA \otimes g"]
      \arrow[rrrr, dashed]
    &&&&
    XA'B
      \arrow[ddr, dashed]
      \arrow[uuullll, phantom, "8.", pos = 0.1, xshift=1.5em]
    \\\\ &
    XAB'
      \arrow[rrrr, swap, "X \otimes f \otimes B'"]
    &&&&
    XA'B'
      \arrow[uulllll, phantom, "\Uparrow X \otimes \otimes_{f,g}", pos = 0.17]
\end{tikzcd}
\]

\bigskip
\[
\begin{array}{llll}
1.\: =\: A \otimes R_{g,X} &
2.\: =\: \otimes_{f,R_{B',X}}  &
3.\: =\:  A' \otimes R_{g,X}   &  
4.\: =\: \otimes_{f,R_{B,X}} \\
5.\: =\: \otimes_{R_{A,X},g}  &
6.\: =\:  R_{f,X} \otimes B'    &
7.\: =\:  \otimes_{R_{A',X},g}  &
8.\: =\:  R_{f,X} \otimes B
\end{array} \]
}
The top cube commutes by $4.(iv),(vi),(viii)$, applied to the 2-morphism 
$A \otimes R_{g,X}$.
The bottom cube commutes by $4.(iv),(vii),(viii)$, applied to the 2-morphism
$R_{f,X} \otimes B$.


We have to verify that these data satisfy the conditions $4.(i) - (viii)$.
These follow from the corresponding conditions holding in $\C$.

\subsection{The Braiding}
$(\bullet \otimes \bullet)$: 
For any two objects we have the morphism
\[(R_{A,B},R_{R_{A,B},\stri}) : 
(A \otimes B ,R_{A,\stri}\otimes R_{B,\stri},\tilR_A \otimes \tilR_B) \to
(B \otimes A,R_{B,\stri} \otimes R_{A,\stri},\tilR_B \otimes \tilR_A)\]
in $\Z(\C)$, where the 2-morphism $R_{R_{A,B},X}$ is defined to be the pasting:

\bigskip 

\[
\begin{tikzcd}[column sep = large, row sep = huge]
    ABX
      \arrow[rr, "R_{A,B} \otimes X"]
      \arrow[d, swap, "A \otimes R_{B,X}"]
      \arrow[drr]
    &&
    BAX
      \arrow[d, "B \otimes R_{A,X}"]
      \arrow[ddll, phantom, "\zwei \tilR_{(A|B,X)}", pos = 0.1, xshift=-1.5em]
      \arrow[ddll, phantom, "\zwei R_{A,R_{B,X}}^{-1}"]
      \arrow[ddll, phantom, "\zwei \tilR_{(A|X,B)}^{-1}", pos = 0.9, xshift=1.5em]
    \\
    AXB
      \arrow[d, swap, "R_{A,X} \otimes B"]
      \arrow[drr]
    &&
    BXA
      \arrow[d, "R_{B,X} \otimes A"]
    \\
    XAB
      \arrow[rr, swap, "X \otimes R_{A,B}"]
    &&
    XBA
\end{tikzcd}
\]

\medskip

First we have to show that $R_{R_{A,B},\stri}$ satisfies $({\to} \otimes {\to})$
and hence is a modification.
This is shown in the following diagram (or follows from the fact that
it is a pasting of modifications).

\[
\begin{tikzcd}[column sep = small]
    ABX
      \arrow[ddd, swap, "A \otimes R_{B,X}"]
      \arrow[rrrr, "R_{A,B} \otimes X"]
      \arrow[dddrrrr, dashed, start anchor = south east, end anchor = north west]
      \arrow[ddr]
    &&& \phantom{ABX} &
    BAX
      \arrow[ddr, "BA \otimes f"]
      \arrow[ddd, dashed]
    \\\\ &
    ABX'
      \arrow[ddd]
      \arrow[rrrr, "R_{A,B} \otimes X"]
      \arrow[dddrrrr, start anchor = south east, end anchor = north west]
    &&&&
    BAX'
      \arrow[ddd, "B \otimes R_{A,X'}"]
    \\
    AXB
      \arrow[ddd, swap, "R_{A,X} \otimes B"]
      \arrow[ddr]
      \arrow[dddrrrr, dashed, start anchor = south east, end anchor = north west]
    &&&&
    BXA
      \arrow[ddd, dashed]
      \arrow[ddr, dashed]
    \\\\ &
    AX'B
      \arrow[ddd]
      \arrow[dddrrrr, start anchor = south east, end anchor = north west]
    &&&&
    BX'A
      \arrow[ddd, "R_{B,X'} \otimes A"]
    \\
    XAB
      \arrow[ddr, swap, "f \otimes AB"]
      \arrow[rrrr, dashed]
    &&&&
    XBA
      \arrow[ddr, dashed]
    \\\\ &
    X'AB
      \arrow[rrrr, swap, "X' \otimes R_{A,B}"]
    &&&&
    X'BA
\end{tikzcd}
\]

\bigskip

The front and the back side of the cube are the 2-morphisms  
$R_{R_{A,B},X}$ and  $R_{R_{A,B},X'}$, respectively.
The top and the bottom are $\otimes_{R_{A,B},f}$ and
$\otimes_{f,R_{A,B}}$, respectively.
The left and the right side are the 2-morphisms corresponding to the 
pseudonatural transformations in the tensor product of the objects $A$ and $B$, 
$(R_A \otimes R_B)_f$ and $(R_B \otimes R_A)_f$, respectively. 

The top triangular prism commutes by 
$(\bullet \otimes (\bullet \otimes {\to}))$.
The bottom triangular prism commutes by
$(\bullet \otimes ({\to} \otimes \bullet))$.
The cube in the middle commutes by 
$(\bullet \otimes (\bullet \otimes {\to}))'$, which is a consequence of
$(\bullet \otimes {\zwei})$ and $(\bullet \otimes {\to} \to)$ as indicated in
Lemma~\ref{Lcub}.

Next, to show that we have really defined a morphism in $\Z(\C)$, we have to verify 
$({\to} \otimes (\bullet \otimes \bullet))$. 
This means we have to check the commutativity of the following diagram.

\vbox{
\[
\begin{tikzcd}[column sep = 0.8 em, row sep = 1.0 em]
    ABXY
      \arrow[drrr, "A \otimes R_{B,XY}"]
      \arrow[ddr]
      \arrow[dddd, swap, "R_{A,B} \otimes XY"]
      \arrow[dddrrrr, phantom, "\Uparrow A \otimes \tilR_{(B|X,Y)}", pos = 0.31]
      \arrow[ddddddddrr, phantom, "1." pos = 0.05]
    \\
    &&&
    AXYB
      \arrow[drrr, "R_{A,XY} \otimes B"]
      \arrow[ddr]
    \\
    &
    AXBY
      \arrow[urr]
      \arrow[ddr, swap, "R_{A,X} \otimes BY"]
      \arrow[drrr, phantom, "\Uparrow \otimes_{R_{A,X},R_{B,Y}}^{-1}"]
    &&&& \phantom{ABXY} &
    XYAB
      \arrow[dddd, "XY \otimes R_{A,B}"]
      \arrow[lllll, phantom, "\Uparrow \tilR_{(A|X,Y)} \otimes B", pos = 0.31]
      \arrow[ddddddllll, phantom, "2.", pos = 0.05]
    \\
    &&&&
    XAYB
      \arrow[urr]
    \\
    BAXY
      \arrow[drrr, dashed]
      \arrow[ddr, swap, "B \otimes R_{A,X} \otimes Y"]    
    &&
    XABY
      \arrow[urr, swap, "XA \otimes R_{B,Y}"]
      \arrow[dddd]
    \\
    &&&
    BXYA
      \arrow[ddr, dashed]
      \arrow[drrr, dashed]
    \\
    &
    BXAY
      \arrow[urr, dashed]
      \arrow[ddr, swap, "R_{B,X} \otimes AY"]
    &&&&&
    XYBA
      \arrow[uuull, phantom, "3.", pos = 0.05, yshift=1.5ex]
    \\
    &&&&
    XBYA
      \arrow[urr, swap, "X \otimes R_{B,Y} \otimes A"]
    \\
    &&
    XBAY
      \arrow[urr, swap, "XB \otimes R_{A,Y}"]
\end{tikzcd}
\]

\bigskip
\[ 1.\: =\: R_{R_{A,B},X} \otimes Y \qquad
2.\: =\: X \otimes  R_{R_{A,B},Y}   \qquad
3.\: =\: R_{R_{A,B},X \otimes Y)}  \]
}

As shown in the diagram below, we decompose this diagram
in the following way:
1) Three tetrahedra of the form
$(\bullet \otimes (\bullet \otimes \bullet \otimes \bullet))$. 
2) One prism of the form
$(\bullet \otimes( \bullet \otimes {\to}))$, namely
$(A \otimes (X \otimes (BY \to YB)))$ (second row, right).
3)  One  prism of  the form
$(\bullet \otimes ({\to} \otimes \bullet))$, namely
$(A \otimes (( BX \to XB) \otimes Y))$ (second row, left).
4)  
One prism of the form $(\bullet \otimes {\zwei})$, namely
$( A \otimes \tilR_{(B | X,Y)})$ (in the middle of the first row). 
All of these diagrams commute by our assumptions.


\noindent
\begin{tikzpicture}
\node (ABXY1) at (1,12)      {$ABXY$};
\node (BAXY1) at (0,9.5)     {$BAXY$};
\node (BXAY1) at (1,7)       {$BXAY$};
\node (BXYA1) at (4,7)       {$BXYA$};

\node (ABXY2) at (2,12.5)    {$ABXY$};
\node (AXYB2) at (6,13)      {$AXYB$};
\node (AXBY2) at (4,11.5)    {$AXBY$};
\node (BXYA2) at (5,7.5)     {$BXYA$};
\node (XYBA2) at (9,8)       {$XYBA$};
\node (XBYA2) at (7,6.5)     {$XBYA$};

\node (AXYB3) at (7.5,12.5)  {$AXYB$};
\node (XYAB3) at (11.5,14)   {$XYAB$};
\node (XAYB3) at (10,11)     {$XAYB$};
\node (XYBA3) at (10.5,7.5)  {$XYBA$};

\node (ABXY4) at (1.5,6)     {$ABXY$};
\node (AXBY4) at (3.5,5)     {$AXBY$};
\node (BXAY4) at (1.5,1)     {$BXAY$};
\node (BXYA4) at (4.5,1)     {$BXYA$};
\node (XBAY4) at (3.5,0)     {$XBAY$};
\node (XBYA4) at (6.5,0)     {$XBYA$};

\node (AXBY5) at (7.3,4.5)   {$AXBY$};
\node (XABY5) at (9.3,3)     {$XABY$};
\node (XBAY5) at (7.3,-0.5)  {$XBAY$};
\node (XBYA5) at (10.3,-0.5) {$XBYA$};

\node (AXBY6) at (8.6,5.5)   {$AXBY$};
\node (AXYB6) at (10.6,7)    {$AXYB$};
\node (XABY6) at (10.6,4)    {$XABY$};
\node (XAYB6) at (12.6,5.5)  {$XAYB$};
\node (XBYA6) at (11.6,0.5)  {$XBYA$};
\node (XYBA6) at (13.6,2)    {$XYBA$};

\draw[->]        (ABXY1) to (BAXY1);
\draw[->]        (ABXY1) to (BXAY1);
\draw[->]        (ABXY1) to (BXYA1);
\draw[->]        (BAXY1) to (BXAY1);
\draw[dashed,->] (BAXY1) to (BXYA1);
\draw[->]        (BXAY1) to (BXYA1);

\draw[->]        (ABXY2) to (BXYA2);
\draw[->]        (ABXY2) to (AXBY2);
\draw[->]        (ABXY2) to (AXYB2);
\draw[->]        (BXYA2) to (XBYA2);
\draw[dashed,->] (BXYA2) to (XYBA2);
\draw[->]        (AXBY2) to (XBYA2);
\draw[->]        (AXBY2) to (AXYB2);
\draw[->]        (XBYA2) to (XYBA2);
\draw[->]        (AXYB2) to (XYBA2);

\draw[->]        (AXYB3) to (XAYB3);
\draw[->]        (AXYB3) to (XYAB3);
\draw[->]        (AXYB3) to (XYBA3);
\draw[->]        (XAYB3) to (XYAB3);
\draw[->]        (XAYB3) to (XYBA3);
\draw[->]        (XYAB3) to (XYBA3);

\draw[->]        (ABXY4) to (AXBY4);
\draw[->]        (ABXY4) to (BXAY4);
\draw[dashed,->] (ABXY4) to (BXYA4);
\draw[->]        (AXBY4) to (XBAY4);
\draw[->]        (AXBY4) to (XBYA4);
\draw[dashed,->] (BXAY4) to (BXYA4);
\draw[->]        (BXAY4) to (XBAY4);
\draw[->]        (XBAY4) to (XBYA4);
\draw[dashed,->] (BXYA4) to (XBYA4);

\draw[->]        (AXBY5) to (XABY5);
\draw[->]        (AXBY5) to (XBAY5);
\draw[dashed,->] (AXBY5) to (XBYA5);
\draw[->]        (XABY5) to (XBAY5);
\draw[->]        (XABY5) to (XBYA5);
\draw[->]        (XBAY5) to (XBYA5);

\draw[->]        (AXBY6) to (AXYB6);
\draw[->]        (AXBY6) to (XABY6);
\draw[->]        (AXBY6) to (XBYA6);
\draw[->]        (AXYB6) to (XAYB6);
\draw[dashed,->] (AXYB6) to (XYBA6);
\draw[->]        (XABY6) to (XAYB6);
\draw[->]        (XABY6) to (XBYA6);
\draw[->]        (XBYA6) to (XYBA6);
\draw[->]        (XAYB6) to (XYBA6);
\end{tikzpicture}

\bigskip

$({\to} \otimes \bullet)$:
For any 1-morphism $(f,R_{f,\stri})\maps (A,R_A,\tilR_A) \to 
(A',R_{A'},\tilR_{A'})$ and any object $(B,R_B,\tilR_B) \in \Z(\C)$ we 
have a 2-isomorphism
 \[ R_{f,B} : (f \otimes B) R_{A',B} \To R_{A,B}(B \otimes f)\]

The following diagram shows that $R_{f,B}$ satisfies 
$({\zwei} \otimes \bullet)$ and is therefore a 2-morphism in $\Z(\C)$.

\[
\begin{tikzcd}[row sep = large, /tikz/column 6/.append style={anchor=base west}]
    ABX
      \arrow[rrr, "f \otimes BX"]
      \arrow[dd, swap, "A \otimes R_{B,X}"]
      \arrow[dr]
      \arrow[dddr]
      \arrow[drrrr, phantom, "\zwei R_{f,B} \otimes X", pos = 0.2]
    &&&
    A'BX
      \arrow[dr, "R_{A',B} \otimes X"]
      \arrow[dd, dashed]
      \arrow[dddr, dashed]
    \\
    &
    BAX
      \arrow[rrr, "B \otimes f \otimes X", pos = 0.42]
      \arrow[dd]
    &&&
    BA'X
      \arrow[dd, "B \otimes R_{A',X}"]
    & ({\to} \otimes (\bullet \otimes \bullet))
    \\
    AXB
      \arrow[dd, swap, "R_{A,X} \otimes B"]
      \arrow[dddr]
      \arrow[rrr, dashed]
    &&&
    A'XB
      \arrow[dd, dashed]
      \arrow[dddr, dashed]
    \\
    &
    BXA
      \arrow[rrr]
      \arrow[dd]
    &&&
    BXA'
      \arrow[dd, "R_{B,X} \otimes A'"]
    & ({\to} \otimes {\to})
    \\
    XAB
      \arrow[dr, swap, "X \otimes R_{A,B}"]
      \arrow[rrr, dashed]
    &&&
    XA'B
      \arrow[dr, dashed]
    \\
    &
    XBA
      \arrow[rrr, swap, "XB \otimes f", pos = 0.42]
    &&&
    XBA'
      \arrow[ullll, phantom, "\zwei X \otimes R_{f,B}", pos = 0.2]
    & ({\to} \otimes (\bullet \otimes \bullet))
\end{tikzcd}
\]
\medskip

The left and right sides are the 2-morphisms $R_{R_{A,B},X}$   and
$R_{R_{A',B},X}$, respectively.  
The front and the back sides are pastings as in our treatment in
Section~\ref{monoidal} of the tensor product of an object and a morphism in $\Z(\C)$.  

We decompose this cube into two commutative triangular prisms of the form 
$({\to} \otimes (\bullet \otimes \bullet))$, correspoding to $(f \otimes (B \otimes X))$
and $(f \otimes (X \otimes B))$, and one cube of the form $({\to} \otimes {\to})$, 
namely $(A \to A' \otimes BX \to XB)$.
\bigskip

$(\bullet \otimes {\to})$:
For any 1-morphism $(g,R_{g,\stri}): (B,R_B,\tilR_B) \to 
(B',R_{B'},\tilR_{B'})$ and any object $(A,R_A,\tilR_A) \in \Z(\C)$, 
we have a 2-iso 
\[ R_{A,g} : (A \otimes g) R_{A,B'} \To R_{A,B}(g \otimes A)\]

The following diagram shows that $R_{A,g}$ satisfies
$({\zwei} \otimes \bullet)$ and is thus a 2-morphism in $\Z(\C)$.

\[
\begin{tikzcd}[row sep = large, /tikz/column 6/.append style={anchor=base west}]
    ABX
      \arrow[rrr, "A \otimes g \otimes X"]
      \arrow[dd, swap, "A \otimes R_{B,X}"]
      \arrow[dr]
      \arrow[dddr]
      \arrow[drrrr, phantom, "\zwei R_{A,g} \otimes X", pos = 0.2]
    &&&
    AB'X
      \arrow[dr, "R_{A,B'} \otimes X"]
      \arrow[dd, dashed]
      \arrow[dddr, dashed]
    \\
    &
    BAX
      \arrow[rrr, "g \otimes AX", pos = 0.42]
      \arrow[dd]
    &&&
    B'AX
      \arrow[dd, "B' \otimes R_{A,X}"]
    & (\bullet \otimes ({\to} \otimes \bullet))
    \\
    AXB
      \arrow[dd, swap, "R_{A,X} \otimes B"]
      \arrow[dddr]
      \arrow[rrr, dashed]
    &&&
    AXB'
      \arrow[dd, dashed]
      \arrow[dddr, dashed]
    \\
    &
    BXA
      \arrow[rrr]
      \arrow[dd]
    &&&
    B'XA
      \arrow[dd, "R_{B',X} \otimes A"]
    & (\bullet \otimes ({\to} \otimes \bullet)')
    \\
    XAB
      \arrow[dr, swap, "XR_{A,B}"]
      \arrow[rrr, dashed]
    &&&
    XAB'
      \arrow[dr, dashed]
    \\
    &
    XBA
      \arrow[rrr, swap, "X \otimes g \otimes A", pos = 0.42]
    &&&
    XB'A
      \arrow[ullll, phantom, "\zwei X \otimes R_{A,g}", pos=0.2]
    & (\bullet \otimes ({\to} \otimes \bullet))
\end{tikzcd}
\]

The decomposition is similar to the one before.
\bigskip

$((\bullet \otimes \bullet) \otimes \bullet)$:
For any objects $(A,R_A,\tilR_A),(B,R_B,\tilR_B),(C,R_C,\tilR_C) 
\in \Z(\C)$ we have the 2-isomorphism
$\tilR_{(A,B|C)}  := 1_{(R_A \otimes R_B)_{C}}$:

\[
\begin{tikzcd}[row sep = huge, column sep = tiny]
    A \otimes B \otimes C
    \arrow[rr, "(R_A \otimes R_B)_C"]
    \arrow[dr, swap, "A \otimes R_{B,C}"]
    & \phantom{1} &
    C \otimes A \otimes B
    \\ &
    A \otimes C \otimes B
    \arrow[u, phantom, "\Uparrow 1", pos = 0.65]
    \arrow[ur, swap, "R_{A,C} \otimes B"]
\end{tikzcd}
\]

\smallskip

$((\bullet \otimes (\bullet \otimes \bullet))$:
For any objects $(A,R_A,\tilR_A),(B,R_B,\tilR_B),(C,R_C,\tilR_C) 
\in \Z(\C)$ we have the 2-isomorphism
$\tilR_{(A|B,C)}$:

\[
\begin{tikzcd}[row sep = huge, column sep = tiny]
    A \otimes B \otimes C
    \arrow[rr, "R_{A, (B \otimes C)}"]
    \arrow[dr, swap, "R_{A,B} \otimes C"]
    & \phantom{1} &
    B \otimes C \otimes A
    \\ &
    B \otimes A \otimes C
    \arrow[u, phantom, "\Uparrow \tilR_{(A|B,C)}", pos = 0.7]
    \arrow[ur, swap, "B \otimes R_{A,C}"]
\end{tikzcd}
\]

To verify that $\tilR_{(A|\stri , \stri)}$ is a 2-morphism in 
$\Z(\C)$, we have to check $( {\zwei} \otimes \bullet)$.
The next diagram gives the proof.

\[
\begin{tikzcd}[row sep = large, column sep = small, /tikz/column 5/.append style={anchor=base west}]
    ABCX
      \arrow[rr, "R_{A,BC} \otimes X"]
      \arrow[dr]
      \arrow[dddrr, dashed, start anchor = {[xshift=0.2em]south}, end anchor = {[xshift=0.5em]north west}]
      \arrow[ddd, swap, "AB \otimes R_{C,X}"]
    &&
    BCAX \arrow[ddd, "BC \otimes R_{A,X}"]
    \\ &
    BACX
      \arrow[ur]
      \arrow[ddr]
      \arrow[ddd]
    &&& (\bullet \otimes (\bullet \otimes \bullet \otimes \bullet))
    \\\\
    ABXC
      \arrow[dr]
      \arrow[ddd, swap, "A \otimes R_{B,X} \otimes C"]
      \arrow[ddddr, end anchor ={[xshift=-0.2em]}]
      \arrow[dddrr, dashed, start anchor = {[xshift=0.2em]south}, end anchor = {[xshift=0.5em]north west}]
    &&
    BCXA \arrow[ddd, "B \otimes R_{C,X} \otimes A"]
    \\ &
    BAXC
      \arrow[ddd]
      \arrow[ddr]
    &&& (\bullet \otimes (\bullet \otimes {\to}) = (A \otimes (B \otimes (CX \to XC)))
    \\\\
    AXBC
      \arrow[ddd, swap, "R_{A,X} \otimes BC"]
      \arrow[ddddr, end anchor ={[xshift=-0.2em]}]
      \arrow[dddrr, dashed, start anchor = {[xshift=0.2em]south}, end anchor = {[xshift=0.5em]north west}]
    &&
    BXCA \arrow[ddd, "R_{B,X} \otimes CA"]
    && (\bullet \otimes (\bullet \otimes \bullet \otimes \bullet))
    \\ &
    BXAC
      \arrow[ur]
      \arrow[ddd]
    &&& (\bullet \otimes ({\to} \otimes \bullet)) = (A \otimes (BX \to XB) \otimes C)
    \\\\
    XABC
      \arrow[dr]
      \arrow[rr, dashed]
    &&
    XBCA
    && (\bullet \otimes (\bullet \otimes \bullet \otimes \bullet))
    \\ &
    XBAC \arrow[ur]
\end{tikzcd}
\]

\bigskip
\bigskip
\medskip

The top triangle corresponds to the 2-morphism 
$\tilR_{(A|B,C)} \otimes X$.
The bottom triangle corresponds to the 2-morphism 
$X \otimes \tilR_{(A|B,C)}$.
The back side is $R_{R_{A,B\otimes C},X}$, the left front side is
$R_{R_{A,B}\otimes C,X}$ and the right front side is
$R_{B \otimes R_{A,C},X}$. 
The decomposition is indicated in the diagram.

\bigskip

Now we have to verify that these data satisfy all the axioms of a braided 
monoidal 2-category.
The tetrahedron
$(\bullet \otimes(\bullet \otimes \bullet \otimes \bullet))$ commutes
by the definition of the objects of $\Z(\C)$.
The diagram
$((\bullet \otimes \bullet) \otimes (\bullet \otimes \bullet))$ 
commutes by the definition of the tensor product of two objects 
in $\Z(\C)$. 
By the same definition can be shown that 
$((\bullet \otimes \bullet \otimes \bullet)\otimes \bullet)$  
commutes in $\Z(\C)$.
Note that because of our special choice of the 2-morphism 
$R_{R_{A,B},\stri}$ that completes the morphism 
$R_{A,B}$ to a morphism in $\Z(\C)$, 
the two 2-morphisms $S^+$ and $S^-$ are equal in $\Z(\C)$.

The other axioms of a braided monoidal 2-category are either part of
our definitions, or else we have indicated within our Remarks which
definitions imply them.  We may summarize by stating:

\begin{thm}  Given any semistrict monoidal category $\C$, the center
$\Z(\C)$ is semistrict braided monoidal 2-category.\end{thm}

%% file: 4embedding.tex
\section{Embedding $\C$ in $\Z(\C)$}

Given a semistrict braided monoidal 2-category $\C$, we would like to
embed it in its center by  
a braided monoidal 2-functor $\F: \C \to \Z(\C)$.  Developing a general
definition of `braided monoidal 2-functor' would require a fair amount of
work.  Luckily, in our case we can restrict ourselves to a very strict
sort of braided monoidal 2-functor which is easy to define.   
The following definition should not be taken as fundamental; it is simply
designed to be the strictest one for which our embedding theorem holds.  

\begin{defn}
Let $(\C,\otimes,I,R,\tilR_{(\stri|\stri ,\stri)}
,\tilR_{(\stri,\stri | \stri)})$ and $(\C',\otimes',I,R',\tilR'_{(\stri|\stri ,\stri)}
,\tilR'_{(\stri,\stri | \stri)})$  be braided monoidal 2-categories. 
A monoidal 2-functor consists of:
\begin{itemize}
\item
A 2-functor $\F : \C \to \C'$
such that $\F(I) = I'$.
\item
A pseudonatural transformation 
\[\xi: (\F \gtimes \F) \circ \otimes' \To \otimes \circ \F , \]
\item an invertible modification
$ \alpha : (1 \otimes \xi)\circ \xi \To (\xi \otimes 1) \circ \xi$,
\end{itemize}
such that the following diagram commutes. \medskip

\vbox{
\[
\begin{tikzcd}[row sep = huge, column sep = tiny] 
    & \F(X)\F(Y)\F(Z)\F(W)
      \arrow [dd, dashed]
      \arrow [rr]
      \arrow [dl]
    &&
    \F(X)\F(YZ)\F(W)
      \arrow [dl]
      \arrow [dd]
        \\
    \F(XY)\F(Z)\F(W)
      \arrow [dd]
      \arrow [rr]
    &&
    \F(XYZ)\F(W)
      \arrow [dd]
      \arrow [dl, phantom, "2.", pos = 0.1]
      \arrow [ul, phantom, xshift=1em, "5.", pos = 0.1]
        \\
    & \F(X)\F(Y)\F(ZW)
      \arrow [dl, dashed]
      \arrow [rr, dashed]
      \arrow [dr, phantom, xshift=-1em, "3.", pos = 0.1]
      \arrow [ur, phantom, "6.", pos = 0.1]
    &&
    \F(X)\F(YZW)
      \arrow [dl]
        \\
    \F(XY)\F(ZW)
      \arrow [rr]
      \arrow [uuur, phantom, "1.", pos = 0.05]
    &&
    \F(XYZW)
      \arrow [uuur, phantom, "4.", pos = 0.05]
\end{tikzcd}
\]

\smallskip

\[ \begin{array}{llll}
1.\hat{=}\: \otimes_{\xi,\xi} & 
2.\hat{=}\:  \alpha_{X\otimes Y,Z,W}& 
3.\hat{=}\:  \alpha_{X,Y,Z\otimes W}  &
4.\hat{=}\: \alpha_{X,Y\otimes Z,W}  \\ 
5.\hat{=}\: \alpha_{X,Y,Z} \otimes \F(W) & 
6.\hat{=}\: \F(X) \otimes \alpha_{Y,Z,W} &
\end{array}    \]
}

\end{defn}

\begin{defn}
A braided monoidal $2$-functor consists of

\begin{itemize}
\item a monoidal 2-functor $(\F,\xi,\alpha)$ 
\item a modification
\[\F_R :\xi \circ \F(R) \To R' \circ \xi, \]
\end{itemize}
such that the following two diagrams commute, expressing the fact that $\F$ respects the modifications 
$\tilR_{(\stri| \stri,\stri)}$ 
and $\tilR_{(\stri,\stri|\stri)}$ up to $\xi$.
\smallskip

\vbox{
\[
\begin{tikzpicture}
\node (xyz)   at (9,10)   {$\F(XYZ)$};
\node (xy-z)  at (4,11.5) {$\F(X)\F(YZ)$};
\node (x-yz)  at (5.3,9)  {$\F(XY)\F(Z)$};
\node (x-y-z) at (0,10)   {$\F(X)\F(Y)\F(Z)$};
\node (xzy)   at (12,6)   {$\F(YXZ)$};
\node (xz-y)  at (7,4.5)  {$\F(Y)\F(XZ)$};
\node (x-zy)  at (7,7.5)  {$\F(YX)\F(Z)$};
\node (x-z-y) at (2.5,6)  {$\F(Y)\F(X)\F(Z)$};
\node (zxy)   at (9,2)    {$\F(YZX)$};
\node (zx-y)  at (4,0)    {$\F(Y)\F(ZX)$};
\node (z-xy)  at (4,3.5)  {$\F(YZ)\F(X)$};
\node (z-x-y) at (0,2)    {$\F(Y)\F(Z)\F(X)$};

\draw[->]        (x-y-z) -> (x-yz);
\draw[->]        (x-y-z) -> (xy-z);
\draw[->]        (x-y-z) -> (x-z-y);
\draw[->]        (x-y-z) -> (z-x-y);
\draw[->]        (x-yz)  -> (xyz);
\draw[->]        (x-yz)  -> (x-zy);
\draw[->]        (xy-z)  -> (xyz);
\draw[dashed,->] (xy-z)  -> (z-xy);
\draw[->]        (xyz)   -> (xzy);
\draw[dashed,->] (xyz)   -> (zxy);
\draw[->]        (x-z-y) -> (x-zy);
\draw[->]        (x-z-y) -> (xz-y);
\draw[->]        (x-z-y) -> (z-x-y);
\draw[->]        (x-zy)  -> (xzy);
\draw[->]        (xz-y)  -> (xzy);
\draw[->]        (xz-y)  -> (zx-y);
\draw[->]        (xzy)   -> (zxy);
\draw[dashed,->] (z-x-y) -> (z-xy);
\draw[->]        (z-x-y) -> (zx-y);
\draw[dashed,->] (z-xy)  -> (zxy);
\draw[->]        (zx-y)  -> (zxy);

\node at (3.5,10.7) {\scriptsize{1.}};
\node at (4.4,11)   {\scriptsize{2.}};
\node at (0.3,9)    {\scriptsize{3.}};
\node at (9.4,8.8)  {\scriptsize{4.}};
\node at (1,9.4)    {\scriptsize{5.}};
\node at (3.7,0.5)  {\scriptsize{6.}};
\node at (4.5,6)    {\scriptsize{7.}};
\node at (6.7,8.8)  {\scriptsize{8.}};
\node at (7.5,3.9)  {\scriptsize{9.}};
\node at (7.6,2)    {\scriptsize{10.}};
\node at (7.5,10)   {\scriptsize{11.}};
\end{tikzpicture}
\]

\smallskip

\[ \begin{array}{llll}
1.\hat{=}\: R'_{(\F(X),\xi)} & 
2.\hat{=}\:  \F_R & 
3.\hat{=}\:  \tilde{R'}_{(\F(X)|\F(Y),\F(Z))}  &
4.\hat{=}\: \F(\tilR_{(X|Y,Z)})  \\ 
5.\hat{=}\: \F_R \otimes \F(Z) & 
6.\hat{=}\: \F(Y) \otimes \F_R &
7.\hat{=}\:  \alpha_{Y,X,Z} &
8.\hat{=}\:  \xi_{R,Z} \\
9.\hat{=}\:  \xi_{Y,R} &
10.\hat{=}\: \alpha_{Y,Z,X}  &
11.\hat{=}\: \alpha_{X,Y,Z} 
\end{array}    \]
}

\vbox{
\[
\begin{tikzpicture}
\node (xyz)   at (9,10)   {$\F(XYZ)$};
\node (xy-z)  at (4,11.5) {$\F(XY)\F(Z)$};
\node (x-yz)  at (5.3,9)  {$\F(X)\F(YZ)$};
\node (x-y-z) at (0,10)   {$\F(X)\F(Y)\F(Z)$};
\node (xzy)   at (12,6)   {$\F(XZY)$};
\node (xz-y)  at (7,4.5)  {$\F(XZ)\F(Y)$};
\node (x-zy)  at (7,7.5)  {$\F(X)\F(ZY)$};
\node (x-z-y) at (2.5,6)  {$\F(X)\F(Z)\F(Y)$};
\node (zxy)   at (9,2)    {$\F(ZXY)$};
\node (zx-y)  at (4,0)    {$\F(ZX)\F(Y)$};
\node (z-xy)  at (4,3.5)  {$\F(Z)\F(XY)$};
\node (z-x-y) at (0,2)    {$\F(Z)\F(X)\F(Y)$};

\draw[->]        (x-y-z) -> (x-yz);
\draw[->]        (x-y-z) -> (xy-z);
\draw[->]        (x-y-z) -> (x-z-y);
\draw[->]        (x-y-z) -> (z-x-y);
\draw[->]        (x-yz)  -> (xyz);
\draw[->]        (x-yz)  -> (x-zy);
\draw[->]        (xy-z)  -> (xyz);
\draw[dashed,->] (xy-z)  -> (z-xy);
\draw[->]        (xyz)   -> (xzy);
\draw[dashed,->] (xyz)   -> (zxy);
\draw[->]        (x-z-y) -> (x-zy);
\draw[->]        (x-z-y) -> (xz-y);
\draw[->]        (x-z-y) -> (z-x-y);
\draw[->]        (x-zy)  -> (xzy);
\draw[->]        (xz-y)  -> (xzy);
\draw[->]        (xz-y)  -> (zx-y);
\draw[->]        (xzy)   -> (zxy);
\draw[dashed,->] (z-x-y) -> (z-xy);
\draw[->]        (z-x-y) -> (zx-y);
\draw[dashed,->] (z-xy)  -> (zxy);
\draw[->]        (zx-y)  -> (zxy);

\node at (3.5,10.7) {\scriptsize{1.}};
\node at (4.4,11)   {\scriptsize{2.}};
\node at (0.3,9)    {\scriptsize{3.}};
\node at (9.4,8.8)  {\scriptsize{4.}};
\node at (1,9.4)    {\scriptsize{5.}};
\node at (3.7,0.5)  {\scriptsize{6.}};
\node at (4.5,6)    {\scriptsize{7.}};
\node at (6.7,8.8)  {\scriptsize{8.}};
\node at (7.5,3.9)  {\scriptsize{9.}};
\node at (7.6,2)    {\scriptsize{10.}};
\node at (7.5,10)   {\scriptsize{11.}};
\end{tikzpicture}
\]

\smallskip

\[ \begin{array}{llll}
1.\hat{=}\: R'_{(\xi,\F(Z))} & 
2.\hat{=}\:  \F_R & 
3.\hat{=}\:  \tilde{R'}_{(\F(A),\F(B)|\F(Z))}  &
4.\hat{=}\: \F(\tilR_{(X,Y|Z)})  \\ 
5.\hat{=}\: \F(X) \otimes \F_R & 
6.\hat{=}\: \F_R \otimes \F(Y) &
7.\hat{=}\:  \alpha_{X,Z,Y} &
8.\hat{=}\:  \xi_{X,R} \\
9.\hat{=}\:  \xi_{R,Y} &
10.\hat{=}\: \alpha_{Z,X,Y}  &
11.\hat{=}\: \alpha_{X,Y,Z} 
\end{array}    \]
}

\end{defn}

\begin{thm} \label{embed} Let $(\C,\otimes,I,T,\tilT_{(\stri|\stri ,\stri)}
,\tilT_{(\stri,\stri | \stri)})$ be a semistrict braided monoidal 2-category, 
and let $\Z(\C)$ be its center.  Then there is a braided monoidal 2-functor 
$\F : \C \to \Z(\C) $ given as follows:
\begin{align*}
    \F(A)      & = (A, T_{A,\stri}, \tilT_{(A|\stri, \stri)}) \\
    \F(f)      & = (f, T_{f,\stri}) \\
    \F(\alpha) & = \alpha
\end{align*}
Moreover $\F$ is injective on objects, morphisms and 2-morphisms, 
and surjective on 2-morphisms.
\end{thm}  

Proof - First let us show that $\F$ is a monoidal 2-functor.  For this, 
we must define a pseudonatural 1-morphism
$\xi_{A,B} : \F(A) \otimes \F(B) \to \F(A \otimes B)$, where 
 $A,B \in \C$. We let
\begin{align*}
    \xi_{A,B} & := (1_{A \otimes B}, \tilT_{(A, B | \stri)}^{-1}) : \\
              & (A, T_{A, \stri}, \tilT_{(A| \stri, \stri)}) \otimes_{\Z(\C)}
                (B, T_{B, \stri}, \tilT_{(B | \stri, \stri)})
                \to
                (A \otimes B, T_{A \otimes B, \stri}, \tilT_{(A \otimes B | \stri, \stri)})
\end{align*}

To be a morphism in $\Z(\C)$,  $\xi_{A,B}$ has to satisfy 
$({\to} \otimes ( \bullet \otimes \bullet))$.
This is equivalent to the axiom 
$((\bullet \otimes \bullet) \otimes (\bullet \otimes \bullet))$ in $\C$.
We show that $\xi$ is natural, not merely pseudonatural.
To this end we first show that for any morphism $f \maps A \to A'$ in $\C$
the following diagram commutes `on the nose'. 
(Remember our shorthand symbol for tensor products in $\Z(\C)$.)
\medskip

\[
\begin{tikzcd}[column sep = 6em, row sep = 10ex]
    {(A \otimes B, T_A \otimes T_B, \tilT_A \otimes \tilT_B)}
    \arrow[r, "{(1_{A \otimes B}, \tilT_{(A,B|\stri)}^{-1})}"]
    \arrow[d, "{(f,T_{f, \stri}) \otimes_{\Z(\C)} B}", swap]
    &
    {(A \otimes B, T_{A \otimes B, \stri}, \tilT_{(A \otimes B|\stri, \stri)})}
    \arrow[d, "{(f \otimes B, T_{f \otimes B, \stri})}"]
    \\
    {(A' \otimes B, T_{A'} \otimes T_B, \tilT_{A'} \otimes \tilT_B)}
    \arrow[r, "{(1_{A' \otimes B}, \tilT_{(A',B|\stri)})}", swap]
    &
    {(A' \otimes B, T_{A' \otimes B, \stri}, \tilT_{(A' \otimes B|\stri, \stri)})}
\end{tikzcd}\]

The morphism ``first right, then down'' equals  
\[ (f \otimes B,T_{(f \otimes B,\stri)} \cdot 
(\tilT^{-1}_{(A, B|\stri)}\circ (\stri \otimes f \otimes B)))\]
The morphism ``first down, then right'' equals
\[ (f \otimes B,((f\otimes B \otimes \stri)\circ \tilT^{-1}_{(A, B|\stri)})
\cdot
(\otimes _{f,T_{(B,X)}}\circ (T_{A',X} \otimes B)) 
\cdot
((A \otimes T_{B,X}) \circ (T_{f,X} \otimes B))) \]
These two $\Z(\C)$-morphisms are equal, since
by $(({\to} \otimes \bullet) \otimes \bullet)$, the underlying 2-morphisms are 
equal:

\medskip

%
%

\[
(({\to} \otimes \bullet) \otimes \bullet)
\qquad
\begin{tikzcd}[row sep = large]
    ABX
    \arrow [rrrr, "T_{AB,X}"]
    \arrow [dd, "f \otimes BX", swap]
    \arrow [drr]
    \arrow [dddrr, phantom, xshift=0.75em, "\Uparrow \otimes_{(f,T_{B,X})}", pos = 0.19]
    \arrow [drrrr, phantom, yshift=1ex, "\Uparrow \tilT_{(A,B|X)}", pos = 0.28]
    & \phantom{5} && \phantom{0} &
    XAB
    \arrow [dd, "X \otimes f \otimes B"]
    \\
    &&
    AXB
    \arrow[urr]
    \arrow[dd]
    \arrow[drr, phantom, xshift=1em, "\Uparrow T_{f,X}B", pos = 0.12]
    &&
    \phantom{3}
    \\
    A'BX
    \arrow[rrrr, dashed]
    \arrow[drr]
    \arrow[drrrr, phantom, yshift=1ex, "\Uparrow \tilT_{(A',B|X)}", pos = 0.28]
    &&&&
    XA'B
    \arrow[uulll, phantom, "\Uparrow T_{f \otimes B,X}", pos = 0.1]
    \\
    &&
    A'XB
    \arrow[urr]
    &&
    \phantom{4}
\end{tikzcd}
\]

\bigskip

Then we must show naturality with respect to morphisms of the form
$g \maps B \to B'$, which is similar.  
Finally, it is easy to show that $\xi$ is also compatible with 2-morphisms.
Using the axiom
$((\bullet \otimes \bullet \otimes \bullet) \otimes \bullet)$
we see that $\xi$ fulfills the associativity condition on the nose, 
so we can define $\alpha$ to be the identity.
\medskip

Next, we show that $\F$ is braided and, in addition, $\F_R = \id$.
 \[ \F_R = \id \maps \xi \circ \F(T) = R^{\Z(\C)} \circ \xi.\]
This is done by the following calculation.
\begin{align*}
    (\xi \circ \F(T) \circ \xi^{-1})_{A,B} & =
    (1_{A \otimes B}, \tilT_{(A,B|\stri)}^{-1}) \circ (T_{A,B},T_{T_{A,B},\stri})
    \circ (1_{A \otimes B}, \tilT_{(A,B|\stri)}) \\
    & = (T_{A,B},{S^-_{A,B,\stri}})              \\
    & = (T_{A,B},{S^+_{A,B,\stri}})              \\
    & = (T_{A,B},R_{(T_{A,B},\stri)})            \\
    & = R^{\Z(\C)}_{A,B}
\end{align*}

Here $\xi^{-1}= (1_{A \otimes B},\tilT_{(A,B|\stri)})$ is the inverse of
$\xi$, as can be easily verified using the composition law for
1-morphisms in $\Z(\C)$.  The third equation holds by our assumption
that $S^+ = S^-$.  The fifth equation holds according to our definition
of the braiding in $\Z(\C)$.

Finally, we must check that both diagrams in the definition of a strong
braided monoidal 2-functor commute.  In the first diagram all the
2-morphisms except $3$ and $4$ are identities.
Note that the 2-morphism labeled $1$, namely $T_{(\F(X),\xi)}$, 
is the identity, since the 1-morphism part of $\xi$ is
the identity, and that by an application of axiom $(\bullet \otimes
{\zwei})$, face $1$ commutes on the nose.  The remaining 2-morphisms $3$
and $4$ are equal.

In the second diagram the 2-morphism $3$ is the identity.  Here, the 2-morphism
$1$ is defined to be $R_{(X,Y|Z)}^{-1}$ and hence agrees with $4$. The 
remaining 2-morphisms are identities, so this diagram also commutes.
\qed

%% file: 5conclusions.tex
\section{Conclusions}

While we have made some progress in understanding 
monoidal 2-categories and braided monoidal 2-categories, it seems
clear that a truly elegant, not to mention correct, treatment of these
concepts requires a better understanding of 3-categories and 4-categories.
In this spirit, we would like to conclude with a list of some issues that
are not yet resolved.

1) We have not included in our definition of semistrict braided monoidal
2-category any axioms involving the unit object $I$ (other than those
appearing in the definition of semistrict monoidal 2-category).  In the
case of a strict braided monoidal category, where $A \otimes I = A
\otimes I = A$ for any object $A$, there are theorems saying that
$R_{A,I} = R_{I,A} = 1_A$.  In the 2-categorical setting the proof for
this theorem turns into an isomorphism: $R_{A,I} \iso R_{I,A} \iso 1_A$.
If we assumed these isomorphisms were equations and in addition that
$R_{f,I} = R_{I,f} = 1_f $ for all 1-morphisms $f : I \-> I$, then we
could conclude that any braided monoidal category with one object gives
a symmetric monoidal category, as expected.

2) In the center $\Z(\C)$ as we have defined it, the 2-morphisms
$\tilR_{(A,B|C)}$ are all identity 2-morphisms, while the 2-morphisms
$\tilR_{(A|B,C)}$ are not.  This points to a curious assymetry in our
definition of the center.  One could equally well have defined the
center so that $\tilR_{(A|B,C)}$ was always the identity, and not
$\tilR_{(A,B|C)}$, but the question is: why is any `symmetry-breaking'
required?  It may be relevant that Gordon, Power and Street's proof
\cite{GPS} that any tricategory is triequivalent to a semistrict
3-category involves a `symmetry-breaking' maneuver.  This occurs because
the definition of semistrict 3-category has an inherent asymmetry, in
that given $f\maps A \to A'$ and $g \maps B \to B'$, the 2-morphism
$\otimes_{f,g}$ goes from $(A\otimes g)(f \otimes B')$ to $(f \otimes
B)(A' \otimes g)$ rather than vice versa.  Perhaps, therefore, the
center construction involves no asymmetries at the level of weak
$n$-categories, but an arbitrary symmetry breaking is needed to
translate it into the framework of semistrict $n$-categories.

Because the strong braided monoidal 2-functor of Theorem~\ref{embed} is
injective on objects, morphisms and 2-morphisms, this result thus serves
as a strictification theorem asserting that any semistrict braided
monoidal 2-category $\C$ is equivalent (in a precise sense) to one for
which $\tilR_{(\stri,\stri|\stri)}$ is trivial.  Indeed, one may prove
this strictification in other ways as well.  One can also, of course,
show that any any semistrict braided monoidal 2-category $\C$ is
equivalent in the same sense to one for which
$\tilR_{(\stri|\stri,\stri)}$ is trivial.